\documentclass[12pt,reqno]{amsart}
\usepackage[utf8]{inputenc}
\usepackage{url}
\usepackage[colorlinks=true, allcolors=blue,linkcolor=blue, citecolor=blue]{hyperref}
\usepackage[capitalise]{cleveref}
\usepackage{amsfonts}
\usepackage{soul}
\usepackage{amsmath}
\usepackage{amssymb}
\usepackage{graphicx}
\newcommand{\abs}[1]{\left|{#1}\right|}
\usepackage{verbatim}
\usepackage{cleveref}
\newcommand{\bbS}{\mathbb{S}}
\newcommand{\norm}[1]{\left\|#1\right\|}
\usepackage{tikz}
\usepackage[shortlabels]{enumitem}
\usepackage{enumitem}
\newcommand{\R}{\mathbb R}
\newlist{Properties}{enumerate}{2}
\setlist[Properties]{label=Property \arabic*., font=\textbf, itemindent=*}
\usepackage{amsthm}

\newtheorem{Theorem}{Theorem}[section]
\newtheorem{Definition}[Theorem]{Definition}
\newtheorem{Lemma}[Theorem]{Lemma}
\newtheorem{Proposition}[Theorem]{Proposition}

\newtheorem{Remark}[Theorem]{Remark}

\newcommand{\intersect}{\cap}
\newcommand{\union}{\cup}
\newcommand{\h}{\mathcal{H}}
\newcommand{\re}{\mathbb{R}}

\newcommand{\mt}{\mathcal{T}}
\newcommand{\mtd}{\mathcal{T}^\delta}

\newcommand{\mf}{\mathcal{F}}

\newcommand{\cone}{\mathcal{C}}
\newcommand{\e}{\varepsilon}
\newcommand{\ms}{\mathcal{S}}

\newcommand{\lt}{\mathcal{L}^3}
\newcommand{\Ln}{\mathcal{L}^n}
\newcommand{\pa}{\partial}

\newcommand{\ssubset}{\subset\joinrel\subset}

\usepackage[mathscr]{euscript}
\DeclareSymbolFont{rsfs}{U}{rsfs}{m}{n}
\DeclareSymbolFontAlphabet{\mathscrsfs}{rsfs}
\numberwithin{equation}{section}

\title{Local minimizers in $3$d of vector Allen-Cahn with a quadruple junction. }
\author{Abhishek Adimurthi}
\address{Department of Mathematics, Indiana University 
Bloomington, IN 47405, USA.}
\email{abadim@iu.edu,abhishek.adimu@gmail.com}
\author{Peter Sternberg}
\address{Department of Mathematics, Indiana University 
Bloomington, IN 47405, USA.}
\email{sternber@iu.edu}

\begin{document}

\begin{abstract}
For $\Omega$ a perturbation of the unit ball in $\R^3$, we establish the existence of a sequence of local minimizers for the vector Allen-Cahn energy. The sequence converges in $L^1$ to a partition of $\Omega$ whose skeleton is given by a tetrahedral cone and thus contains a quadruple point. This is accomplished by proving that the partition is an isolated local minimizer of a weighted perimeter problem arising as the associated $\Gamma$-limit of the sequence of Allen-Cahn functionals. 
\end{abstract}

\maketitle
\pagestyle{myheadings}
\markright{4-phase quadruple junction.}

\section{Introduction.}
For a domain $\Omega\subset\R^3$ described through a slight deformation of a ball, we establish the existence of a sequence of local minimizers $u_\e:\Omega\to\R^3$ to a vector Allen-Cahn (or Modica-Mortola) type energy that converges as $\e\to 0$ to a partition of $\Omega$ exhibiting a quadruple junction. Invoking the well-known machinery of $\Gamma$-convergence for Modica-Mortola to be found in \cite{Baldo,KohnPS}, our approach consists of proving that the partition of $\Omega$ into four subdomains via a tetrahedral cone represents an isolated local minimizer of perimeter in the $L^1$ topology.

To be more specific, we define the sequence of energies
\begin{equation*}
 E_{\varepsilon}(u) := \int_{\Omega} \frac{\varepsilon}{2} |\nabla u|^2 + \frac{1}{\varepsilon} W(u) \,dx, 
\end{equation*}
where $u:\Omega\to\R^3$ and $W : \re^3 \mapsto \re$ is a non-negative function that vanishes on a set of four distinct points $P:= \{ p_1, p_2, p_3, p_4\} \subset \re^3$. For say $C^{1,\alpha}$ potentials $W$, a local minimizer represents a solution to the problem
\begin{equation}\label{1}
    \e^2\Delta u = W_u(u)\;\mbox{in}\;\Omega,\quad\nabla u\cdot\nu_{\pa\Omega}=0\;\mbox{on}\;\pa\Omega,
\end{equation}
where $\nu_{\pa\Omega}$ denotes the outer unit normal to $\pa\Omega$. In the analogous situation in one lower dimension where $u:\Omega\to\R^2$ with $\Omega\subset\R^2$ taken as a certain deformation of a disc, and with $W$ given by a triple-well potential, this program is carried out in \cite{Zeimer} to produce solutions to \eqref{1} possessing triple junction structure.

We note that entire solutions $u:\R^3\to\R^3$ to this PDE exhibiting quadruple junction structure have been constructed under symmetry assumptions on the potential $W$ in \cite{Alikakos1,AlikFusco,GuiSch}. See \cite{Alikakos_Book} for a definitive exposition on entire solutions for symmetric potentials in all dimensions. One motivation for the present article, in work presently being pursued, is to construct an entire, locally minimizing solution whose `blow-downs' converge to a quadruple junction partition of $\R^3$. {\color{black} Here `locally minimizing' is meant in the De Giorgi sense, that is, an entire solution that minimizes the energy on compact sets when compared with competitors sharing its boundary values. Such a strategy for non-symmetric triple-well potentials in two dimensions has been carried out in \cite{AliGeng1,Geng,Sandier} to produce entire solutions $u:\R^2\to\R^2$ with triple junction structure.  The argument in two dimensions proceeds by first taking an appropriate blow-up of a local minimizer to $E_\e$ solving \eqref{1} or else subject to a well-chosen sequence of Dirichlet conditions. The limit of these blow-ups yields an entire, locally minimizing solution to $\Delta u=W_u(u)$. Next one argues that blow-downs from infinity of this entire solution necessarily converge to a locally minimizing partition of the plane. Finally, one shows that this locally minimizing partition consists of three rays emanating from a common point at prescribed angles, that is, a triple junction.}

If a similar strategy were successful in three dimensions, the advantage would be two-fold: it would require no symmetry assumption on $W$ and it would lead to a locally minimizing entire solution in the sense of De Giorgi. The previously constructed entire solutions in the case $u:\R^3\to\R^3$ mentioned above do not come with any guarantee of stability.

This entire program represents an attempt at `desingularization' of all possible minimal structures. From this perspective, the most basic case is a plane whose diffuse counterpart for $\Delta u=W_u(u)$ is a heteroclinic connection, that is, an ODE joining two wells. Next comes the aforementioned triple junction, and here we contemplate the other possible limiting cone, the quadruple junction in $\R^3$, as per Taylor's seminal classification in \cite{Taylor}. {\color{black} In \cite{Taylor}, possible tangent cones are classified as in \cref{Taylor} when all interfaces are weighted equally (that is, all $c_{ij}$'s are equal in \eqref{zz} below). }
\begin{figure}
    \centering
    \includegraphics[width=1.1\linewidth]{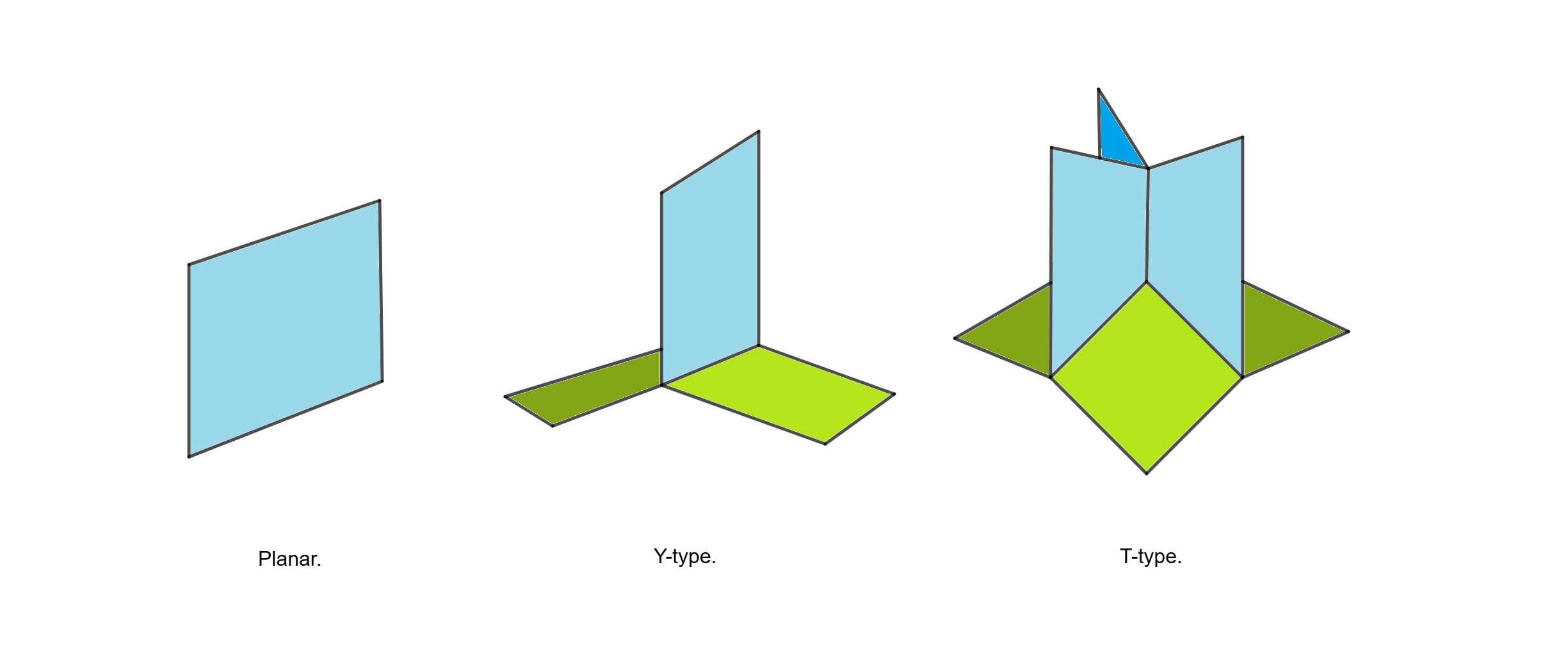}
    \caption{Tangent cones classified in \cite{Taylor} when weights are equal.}
    \label{Taylor}
\end{figure}

 Before describing the partitioning problem in more detail we recall the relevant work on $\Gamma$-convergence.
Generalizing the approach in \cite{FT,PS1,PeterVect} on $\Gamma$-convergence for vector Modica-Mortola with double well-type potentials, under some mild assumptions on a multi-well potential $W$, the results in \cite{Baldo} establish the $\Gamma$-convergence of $E_{\varepsilon}$ in the $L^1$ topology to a functional $E_0$ given in our setting by
\begin{equation}
  E_{0} (v) :=   \sum_{1\leq i<j \leq 4} c_{ij} \h^2 \left( \pa^* R_i \intersect \pa^* R_j  \intersect \Omega \right) \label{zz}
\end{equation}

for $v \in BV(\re^3;P)$. Here the region $R_i$ is defined as
\[
R_i := \{ x \in \Omega \ | \ v(x) = p_i\},
\]
so that $\mathcal{R}:= (R_1,R_2,R_3,R_4)$ represents a partition of $\Omega$ by four sets of finite perimeter. In what follows, we will frequently write $E_0(\mathcal{R})$ rather than $E_0(v)$.
Here $\h^2$ is two-dimensional Hausdorff measure and $\pa^*A$ refers to the reduced boundary of the set $A$. The coefficients $c_{ij}$ are given by 
\begin{equation}
  c_{ij} : = d(p_i, p_j),\label{cij}  
\end{equation} with $d(p,q)$ denoting the metric defined through
\begin{equation*}
    d(p,q) := \inf \left\{ \sqrt{2} \int_0^1 \sqrt{W(\gamma(t))} | \gamma'(t)| dt \left| \begin{aligned}
        &\gamma\in C^1([0,1]; \re^3), \\
        &\gamma(0) = p, \\
        &\gamma(1) = q,
    \end{aligned}\right. \right\}
\end{equation*}
for any $p,q\in\R^3$.

To describe the main result of this paper, Theorem \ref{isolated_local_minimizer}, we consider a tetrahedral cone $\cone$ centered at the origin. The {\it regular} tetrahedral cone, which corresponds to the case where all $c_{ij}$'s in \eqref{cij} are equal, is defined as the cone over a regular tetrahedron. Such a cone partitions $\R^3$ into four regions, say $\tilde{S}_1,\ldots,\tilde{S}_4$, each having boundary consisting of three planes. We will denote the (constant) unit normal to the plane separating region $\tilde{S_i}$ from $\tilde{S_j}$, pointing from $\tilde{S}_i$ to $\tilde{S}_j$, by $n_{ij}$. 

The only point in common to the boundary of all four regions is the origin, but given any triple of distinct indices $\{i,j,k\}\subset\{1,2,3,4\}$, one has that $\pa \tilde{S}_i\cap \pa\tilde{S}_j\cap\pa\tilde{S}_k$ consists of a ray emanating from the origin.
Stationarity of the tetrahedral cone with respect to surface area implies a `balance law' for the normals to the planes meeting at such a ray:
\begin{equation} n_{ij}+n_{jk}+n_{ki}=0\quad\mbox{for every distinct triple}\;i,j,k\in\{1,2,3,4\}.\label{balance}
\end{equation}

More generally, in this article, we will consider partitions arising as critical points of $E_0$ with {\it unequal} weights $c_{ij}$ and the partitions of interest here are induced by {\it non-symmetric tetrahedral cones} centered at the origin. Their description is similar to that of the regular tetrahedral cone, but now stationarity for $E_0$ leads to the requirement that the angles at which three planes meet along the triple-junction like rays are determined through the conditions
\begin{equation}
c_{ij}n_{ij}+c_{jk}n_{jk}+c_{ki}n_{ki}=0\quad\mbox{for every distinct triple}\;i,j,k\in\{1,2,3,4\}.\label{gencij}
\end{equation}

Following \cite[Theorem 1]{Schoenberg_Annals} and \cite{Morgan}, such a non-symmetric cone can be constructed provided the matrix

\begin{equation}\label{cij_matrix}
M:=    \begin{bmatrix}
2c_{12}^2 & c_{12}^2 + c_{13}^2 - c_{23}^2& c_{12}^2 + c_{14}^2 - c_{24}^2 \\
 c_{12}^2 + c_{13}^2 - c_{23}^2 & 2c_{13}^2 & c_{13}^2 + c_{14}^2 - c_{34}^2 \\
 c_{12}^2 + c_{14}^2 - c_{24}^2 & c_{13}^2 + c_{14}^2 - c_{34}^2 & 2c_{14}^2
\end{bmatrix}
\end{equation} 
is positive definite. Under this assumption, there exists a tetrahedron with vertices $A_1,A_2,A_3$ and $A_4$ lying in $\re^3$ with the property that the length of the edge joining $A_i$ and $A_j$ is $c_{ij}$. Then for all distinct pairs $i,j \in \{1,2,3,4\}$, one defines a unit normal vector to the plane separating phase $i$ from phase $j$ via
 \begin{equation}\label{defn_of_normal_to_cone}
 n_{ij} := \frac{1}{c_{ij}}( A_j- A_i).
\end{equation}
 The six planes corresponding to each of these normals $\{n_{ij}\}$ passing through {\color{black} the circumcenter of the tetrahedron $A_1 A_2 A_3 A_4$} lead to the non-symmetric version of the cone $\cone$ satisfying \eqref{gencij}. {\color{black} We recall that for a tetrahedron, the circumcenter is the center of the unique sphere passing through the four vertices. Equivalently, if one considers the plane passing through the midpoint of each edge orthogonally, then the circumcenter is the intersection point of these six planes. In what follows, we will always take the circumcenter as the origin of our coordinate system in $\re^3$. } We also note that positive definiteness of $M$ implies strict triangle inequalities between any three of the $c_{ij}$'s.

{\color{black}At times, it will be more convenient to write the weighted perimeter $E_0(\mathcal{R})$  as
\begin{equation}
E_0(\mathcal{R})=\sum_{i=1}^4 c_i\h^2(\pa^*R_i\cap\Omega), 
\label{ciE}
\end{equation}
where the coefficients $c_i$ are related to the coefficients $c_{ij}$ in \eqref{zz} via the system of linear equations
\begin{equation*}
    c_i+c_j=c_{ij}\quad\text{for}\;1\leq i,j\leq 4,\;i\not=j.
\end{equation*}}
In fact, we will make the stronger assumption that the weights $c_{ij}$ are not too far apart from each other, namely:
\begin{equation}
    c_i>\frac{1}{2}c_j\quad\text{for all distinct}\;i,j\in\{1,2,3,4\}.\label{ciclose}
\end{equation}
 One can verify, e.g. using Mathematica, that condition \eqref{ciclose} implies the positive definiteness condition \eqref{cij_matrix}.

In Section \ref{defn_of_domain} we will describe a domain $\Omega\subset\R^3$ arising as a perturbation of a ball such that the partition $\mathcal{S}:=\{S_1,S_2,S_3,S_4\}$ defined through $S_i:=\tilde{S}_i\cap\Omega$ results in an isolated $L^1$ local minimizer of $E_0$ in the sense that for some $\delta=\delta(\Omega)>0$, the function 
\begin{equation}
  u_0(x):=\sum_{i=1}^4 p_i\chi_{S_i}(x)\label{uzero}  
\end{equation}
 satisfies the condition
\[
E_0(u_0)<E_0(v)\quad\text{provided}\;v\in BV(\Omega;P)\;\text{and}\;0<\norm{v-u_0}_{L^1(\Omega)}<\delta.
\]
This is the content of our main result, Theorem \ref{isolated_local_minimizer}. 
As we shall see, the crucial property of this deformation of the ball is that in a neighborhood of $\cone\cap\pa\Omega$, the vectors $n_{ij}$ satisfy conditions \eqref{goodnormals}.

 In light of Theorem 4.1 of \cite{KohnPS}, Theorem \ref{isolated_local_minimizer} immediately implies:
\begin{Theorem}\label{pde}
    There exists a domain $\Omega\subset\R^3$ such that for all $\e$ sufficiently small, the energy $E_\e$ possesses an $L^1$ local minimizer $u_\e$. As $\e\to 0$, one has that $u_\e\to u_0$ in $L^1(\Omega)$ with $u_0$ given by \eqref{uzero}.
\end{Theorem}

\noindent\underline{Organization of the paper:} We begin in Section \ref{finper}
with some general comments and notation regarding sets of finite perimeter. Section \ref{tools} contains a few tools we will invoke, and in Section \ref{defn_of_domain} we describe the deformation of the ball required for the proof of Theorem \ref{isolated_local_minimizer}. 
The proof of Theorem \ref{isolated_local_minimizer} follows in two steps. In the first, phrased separately as Theorem \ref{lemma1}, we argue that one need only consider partitions that meet $\pa\Omega$ uniformly close to where the partition $\mathcal{S}$ meets $\pa\Omega$. This is accomplished by proving a boundary version of the infiltration lemma (see \cite{Leonardi}, \cite[Lemma 30.2]{Maggi}). It is here that our proof requires the assumption \eqref{ciclose}. The second step then consists of a callibration argument in the spirit of \cite{Morgan} and \cite{Zeimer}. In Section \ref{Construct} we provide an explicit description of a domain supporting an isolated local minimizer $u_0$ as in \eqref{uzero}.

\section{Preliminaries.}
\subsection{Partitions with sets of finite perimeter.}\label{finper}
We recall that a set $E\subset\re^n$ has finite perimeter in an open set $\Omega\subset \re^n$ if the characteristic function $\chi_E$ has bounded variation, and we denote the perimeter within $\Omega$ by ${\rm Per}(E;\Omega)$. We also recall that sets of finite perimeter 
possess a measure-theoretic exterior normal which is suitably general
to ensure the validity of the Gauss-Green theorem. A unit vector
$\nu^*$ is defined 
as the measure theoretic exterior normal to $E$ at $x$ provided
\[
\lim_{r\to0}\frac{\lt\left(B(x,r)\cap\{y:(y-x)\cdot\nu^*<0,y\notin E\}\right)}{r^n}=0
\]
and
\[
\lim_{r\to0}\frac{\lt\left(B(x,r)\cap\{y:(y-x)\cdot\nu^*>0,y\in E\}\right)}{r^n}=0,
\]
where $B(x,r)$ denotes the open ball of radius $r$ centered at $x$ and $\lt$ denotes three-dimensional Lebesgue measure. We denote the set of such points $x$ by $\pa^* E$, the reduced boundary of $E$.

By definition, sets of finite perimeter are determined only up to
sets of measure zero. In order to avoid this ambiguity, whenever a
set of finite perimeter, $E$, is considered, we shall always employ
the measure-theoretic closure as the set to represent $E$. Thus, with
this convention, we have
\[
x\in E\;\hbox{if and only if}\;\limsup_{r\to 0}\frac{\mathcal{L}^n( E\intersect B(x,r))}{\mathcal{L}^n(B(x,r))}>0.
\]

We also note that 
\[
{\rm Per}(E;\Omega)=H^{n-1}(\Omega\cap\partial^{*}E).
\]

Throughout, we will employ the concept of a partition of $\Omega$ by
$4$ sets of finite perimeter.

\begin{Definition}\label{partition_definition}
    We say $\mt :=(T_1, T_2,T_3, T_4)$ is a partition of a bounded domain  $\Omega \subset \re^3$ if each $T_i\subset\Omega$ is a set of finite perimeter, $\mathcal{L}^3 \left( \Omega \setminus \union_{i=1}^3 T_i\right)$\\$=0$, and
   for each distinct pair $i,j\in\{ 1,2,3,4\} $, one has $\mathcal{L}^3(T_i \intersect T_j ) =0$.
\end{Definition}
We recall from the introduction that for a partition $\mathcal{R}= (R_1, R_2,R_3, R_4)$ of $\Omega$, we define
\begin{equation}\label{originalE0}
E_0(\mathcal{R}) := \sum_{1 \leq i < j \leq 4}c_{ij}\h^{2} \left(\pa^* R_i \intersect \pa^* R_j \intersect \Omega\right).
\end{equation}

Further, if $\mathcal{R}':= (R_1',R_2',R_3', R_4')$ is another partition of $\Omega$, then we set
\[
\mathcal{L}^3(\mathcal{R}\Delta \mathcal{R}') := \sum_{i=1}^4 \mathcal{L}^n\Big[(R_i \setminus R_i') \union (R_i' \setminus R_i)\Big].
\]

For more background on sets of finite perimeter we refer the reader to e.g. \cite{EvansGareipy,G,Maggi}.
\subsection{Two helpful tools}\label{tools}
In the analysis to follow, we require the following two results.
\begin{Lemma}\label{hyperplane-theorem}
{\rm (cf. \cite[Lemma 2.2]{Zeimer})}
    Let $E \subset \re^n, n\geq 2$ be a set of finite perimeter. Suppose there is a point $x_0$ in $\pa^* E$ with the property that for some cube $Q_0$ centered at $x_0$, the unit outward normal to $E$ given by $\nu_E(x)$ is a constant $\nu_0$ for all $x \in Q_0 \intersect \pa^* E$. Then, $Q_0 \intersect \pa^* E = Q_0 \intersect \Pi$, where $\Pi$ is the unique hyperplane in $\re^n$ containing $x_0$ with normal $\nu_0$.
\end{Lemma}

\begin{Proposition}[Relative isoperimetric inequality]\label{isoperimetric_theorem}
Let $\Omega \subset \re^n$, $n\geq 2$ be a bounded and connected domain with $C^1$ boundary. 
Then, there exist constants $r_0(\Omega)>0$ and $C_0=C_0(n,\Omega)>0$ such that for any $x_0 \in \pa \Omega$ and any $E \subset \Omega \intersect B(x_0,r_0)$ of finite perimeter, one has:
\begin{equation*}
    {\rm Per}(E;\Omega) \geq C_0 \Big[\mathcal{L}^n (E  )\Big]^{\frac{n-1}{n}}.
\end{equation*}
\end{Proposition} 

A version of \cref{isoperimetric_theorem} for the case when $\Omega$ is a ball in $\re^n$ can be found e.g. in \cite[Section 5.6]{EvansGareipy}. For completeness, we provide a proof of \cref{isoperimetric_theorem} below, though the proof is essentially unchanged.

\begin{proof}  We fix any $r_0>0$ depending on $\Omega$ such that for all $x_0\in\pa \Omega$ one has
\begin{equation}
    \Ln(B(x_0,r_0)\intersect \Omega)\leq \Ln (\Omega\setminus B(x_0,r_0)).\label{goodr}
\end{equation}
For a function $f:\re^n \mapsto \re$, we set \[
    (f)_\Omega := \frac{1}{\mathcal{L}^n(\Omega)} \int_\Omega f(y) dy.
    \]
For any $x_0 \in \pa \Omega$ and any subset $E\subset B(x_0,r_0) \intersect \Omega$, one has
    \begin{equation*}
        \begin{aligned}
            &\int_\Omega \left| \chi_E(y) - (\chi_E)_\Omega \right|^{\frac{n}{n-1}}dy \\
            &= \int_\Omega \left| \chi_E(y) - \frac{\mathcal{L}^n(E)}{\mathcal{L}^n(\Omega)}\right|^{\frac{n}{n-1}}dy\\
            &=\left(\frac{\Ln(\Omega-E)}{\Ln(\Omega)}\right)^{\frac{n}{n-1}} \Ln( E) + \left(\frac{\Ln( E)}{\Ln(\Omega)}\right)^{\frac{n}{n-1}}\Ln(\Omega-E).\\
            & \geq \left(\frac{1}{2}\right)^{\frac{n}{n-1}}\Ln(E),
        \end{aligned}
    \end{equation*}
where the final inequality follows from \eqref{goodr}. Hence,
\begin{equation}
    \lVert  \chi_E - (\chi_E)_\Omega \rVert_{L^{\frac{n}{n-1}}(\Omega)} > \frac{1}{2}\left[ \Ln (E)\right]^{\frac{n-1}{n}}.\label{easy1}
\end{equation}

Now by the Poincar\'e inequality applied to any $f\in W^{1,\frac{n}{n-1}}(\Omega)$, we have
\begin{equation*}
    \left\lVert f -(f)_\Omega \right\rVert_{L^{\frac{n}{n-1}}(\Omega)} \leq C_1 \lVert Df \rVert_{L^{\frac{n}{n-1}}(\Omega)},
\end{equation*}
for some $C_1=C_1(\Omega)$.
Furthermore, we recall that for any $f \in BV(\Omega)$, there exist functions $\{f_k\}_{k=1}^{\infty}\subset BV(\Omega)\intersect C^{\infty}(\Omega)$ such that 
\[
f_k \rightarrow f\;\mbox{ in}\; L^1(\Omega)\quad\mbox{and}\quad 
       \lVert Df_k \rVert(\Omega) \rightarrow \lVert Df\rVert (\Omega)\;\mbox{as}\; k\rightarrow\infty,
   \]
   (see e.g. \cite[Chap. 1]{G}). 
Approximating $f=\chi_E$ by smooth functions we may thus extend the Poincar\'e inequality to such an $f$ with the right-hand side now denoting total variation (perimeter). Combining this inequality with \eqref{easy1}, we complete the proof with $C_0=\frac{1}{2C_1}$.
\end{proof}

\subsection{Description of the Domain.}\label{defn_of_domain}

\begin{figure}
    \centering
    \includegraphics[width=0.9\linewidth]{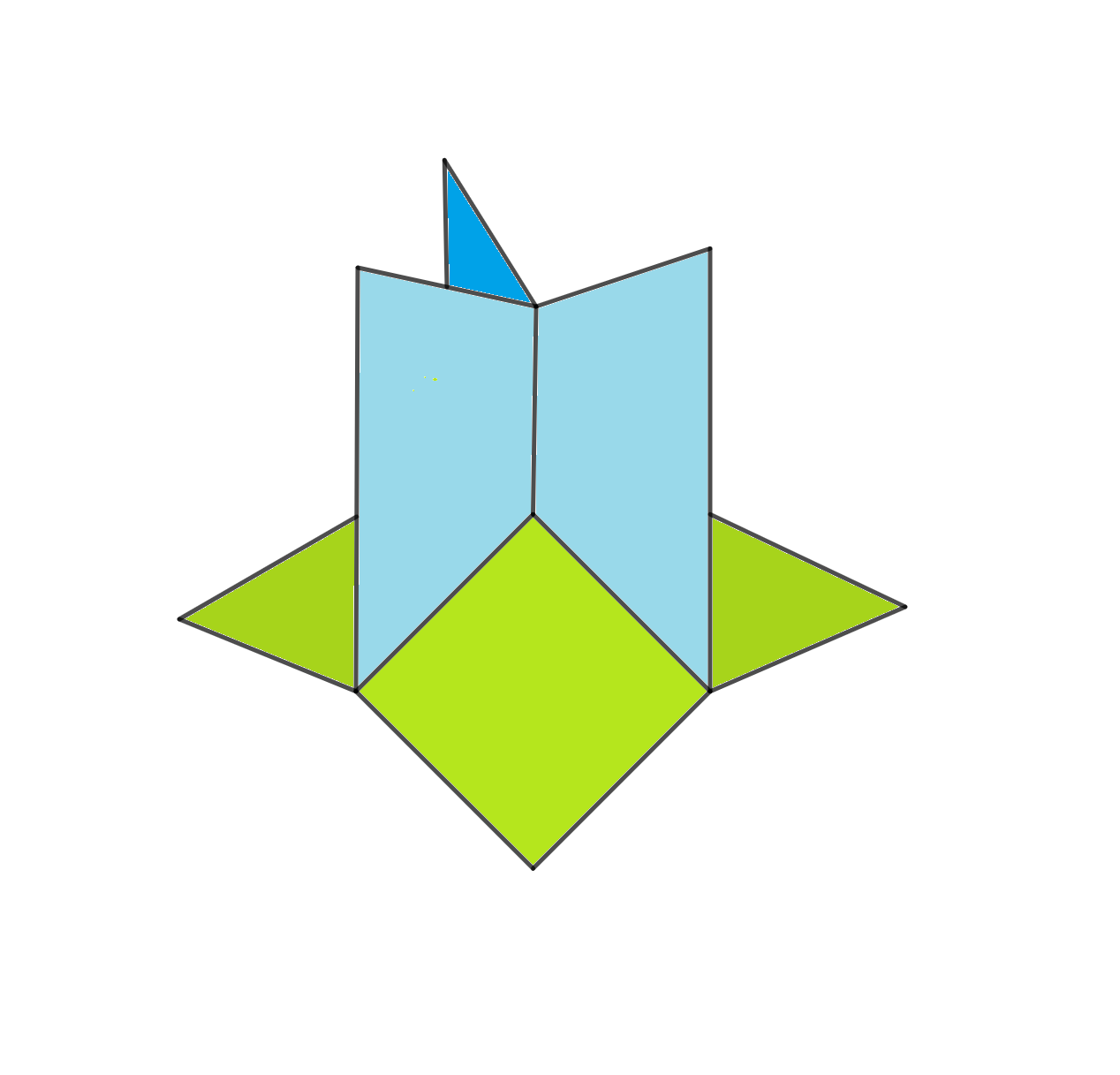}
    \caption{A tetrahedral cone.}
    \label{Tetrahedral_cone}
\end{figure}

To describe the domain $\Omega\subset\R^3$ in which we will construct a sequence of local minimizers of the vector Allen-Cahn system, the key object is the infinite tetrahedral cone ({\color{black} see \cref{Tetrahedral_cone}}).
 We recall that this cone  partitions $\R^3$ about the origin into four regions, say $\tilde{S}_1,\ldots,\tilde{S}_4$, each having boundary consisting of three planes with the normals satisfying the conditions \eqref{balance} in the symmetric case or \eqref{gencij} in the non-symmetric case.

\begin{figure}
	\centering
	\includegraphics[width = 0.6\textwidth]{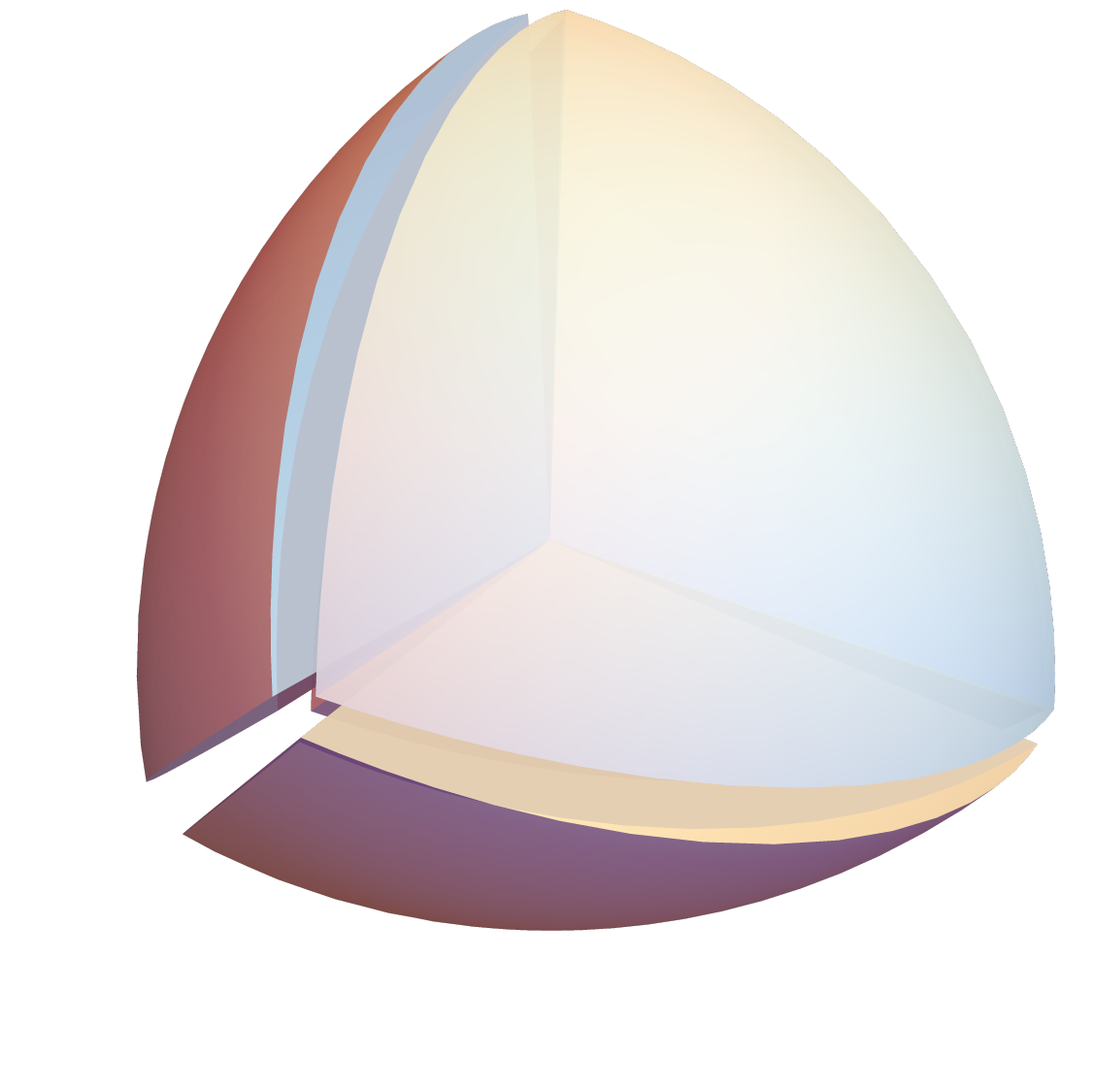}
	\noindent\caption{Partition of the unit ball by the tetrahedral cone. Each connected region here represents one of the sets $\tilde{S}_i$ restricted to the unit ball.}
	\label{Peter1}
\end{figure}
 \begin{figure}
	\centering
	\includegraphics[width = 0.6\textwidth]{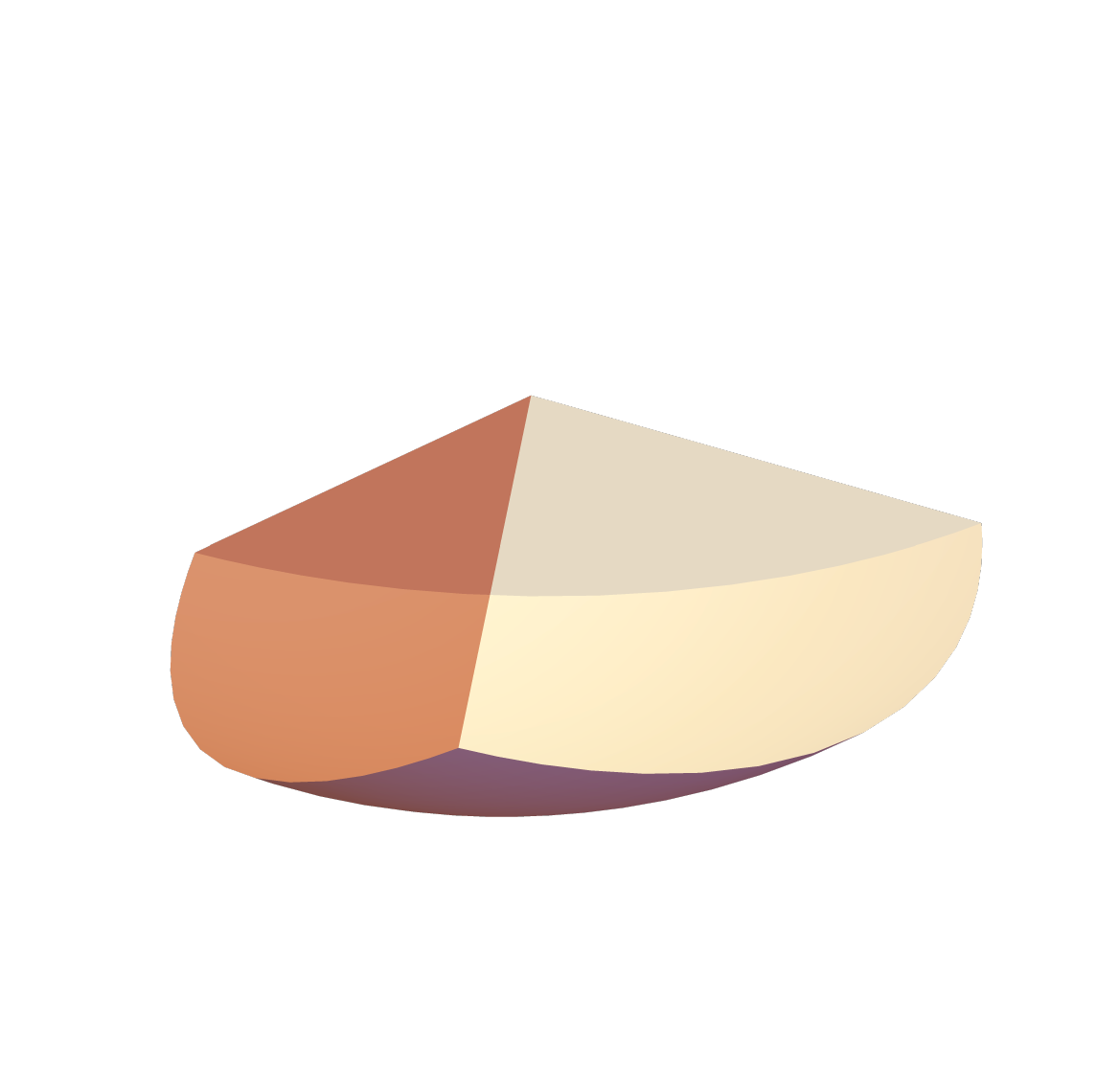}
	\caption{One of the components of the partition of the unit ball depicted in \cref{Peter1}.}
	\label{Peter2}
\end{figure}

Let us now focus on the restriction of this $\cone$ to the unit ball in $\R^3$ centered at the origin.
For each distinct triple of indices, we denote by $\tilde{q}_{ijk}$ the point where the ray on $\cone$ emanating from the origin lying on $\pa\tilde{S}_i\cap\pa\tilde{S}_j\cap\pa\tilde{S}_k$ meets $\bbS^2$. The intersection of $\cone$ with $\bbS^2$ consists of four geodesic triangles, having three of the four $\tilde{q}_{ijk}$'s as vertices.  Furthermore, there are a total of six circular arcs, namely the sides of these triangles, comprising $\cone\cap\bbS^2$, and we denote these arcs by $\tilde{\gamma}_{ij}$, so that $\tilde{\gamma}_{ij}$ has endpoints $\tilde{q}_{ijk}$ and $\tilde{q}_{ijl}$ for distinct indices $i,j,k$ and $l$ with
\[
\tilde{\gamma}_{ij}=\pa \tilde{S}_i\cap\pa \tilde{S}_j\cap \bbS^2.
\]

\begin{figure}
    \centering
    \includegraphics[width=0.5\linewidth]{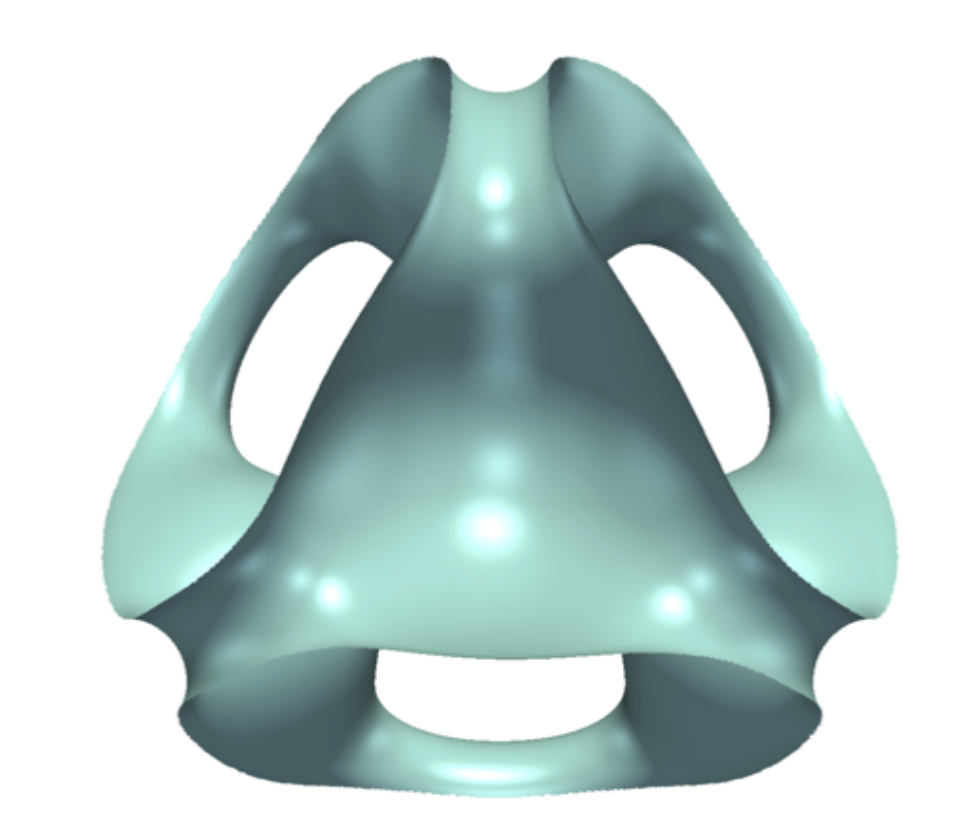}
    \caption{The `troughs' and `valleys' on the surface of the deformed sphere.}
    \label{POI}
\end{figure}

The restriction of $\{\tilde{S}_1,\ldots,\tilde{S}_4\}$ to the unit ball constitutes a partition of the unit ball, see Figures \ref{Peter1} and \ref{Peter2}. We now describe a $C^1$ deformation of $\bbS^2$ in a neighborhood of 
\[  
\cone\cap \bbS^2=\bigcup_{\substack{i,j\in\{1,2,3,4\}\\i\not=j}}\tilde{\gamma}_{ij}.
\]
The resulting surface forms the boundary of a domain $\Omega$ which is partitioned by $\cone$ to form an isolated local minimizer of $E_0$ in the domain.
This deformation is carried out so as to create six troughs along neighborhoods of each $\tilde{\gamma}_{ij}$, except near its endpoints $\tilde{q}_{ijk}$ and $\tilde{q}_{ijl}$. By `troughs' we mean that within this neighborhood of each $\tilde{\gamma}_{ij}$, the closest points to the origin are those on the circular arc itself, with the deformed surface sloping downwards towards $\tilde{\gamma}_{ij}$ on either side. Then, in a neighborhood of each of the four triple junctions $\tilde{q}_{ijk}$ where three such troughs meet, we depress $\bbS^2$ to create small valleys, for example hemispheres. The deepest point of each valley we denote by $q_{ijk}$ where, say, $q_{ijk}=(1-\lambda)\,\tilde{q}_{ijk}$ for some small $\lambda>0$. See \cref{POI}. Finally, away from these valleys and troughs we glue the surface in a $C^1$ manner to a slightly enlarged sphere, say of radius $1+\lambda$, though as we will see, the particular shape away from neighborhoods of the troughs and valleys is unimportant.
Of course, such a construction is not unique, but we describe in detail one such construction in Section \ref{Construct}. 

This deformation of $\bbS^2$ now bounds an open set that we denote by $\Omega$. The tetrahedral cone $\cone$ partitions $\Omega$ into four subsets that we henceforth denote by $S_1,\ldots, S_4$, where for each distinct $i$ and $j$, $\pa S_i\cap \pa S_j$ and $\pa \tilde{S}_i\cap\pa\tilde{S}_j$ lie on the same plane.
 Finally, we let
\[
\gamma_{ij}:=\pa S_i\cap\pa S_j\cap \pa\Omega.
\]

Ultimately, what we require of our domain $\Omega$ is that $\cone$ determines a partition $\Omega$ given by $\mathcal{S}:=\{S_1,S_2,S_3,S_4\}$ satisfying the following property:
\vskip.1in
There exists a real number $\eta=\eta(\Omega)>0$ such that for each distinct pair of indices $i,j\in\{1,2,3,4\}$ and every point in the set
\[
     \{ x \in \pa S_i \intersect \pa \Omega :\, {\rm dist} ( x, \gamma_{ij}) < \eta\},
    \]
one has the property:
\begin{equation}
    \left\{\begin{aligned}\label{goodnormals}
        &n_{ij}\cdot\nu_{\pa\Omega}(x) >0\quad\text{if}\;x\not\in \gamma_{ij},\\
        &n_{ij}\cdot \nu_{\pa\Omega}(x) = 0 \quad\text{if}\; x\in\gamma_{ij} .
    \end{aligned}\right.
    \end{equation}
{\color{black} We note that as the cones become more non-symmetric, one must consider smaller $\eta$ to ensure that these valleys are all disjoint, as are the various troughs.}

\begin{Remark}
    We note that the arguments to follow only require that the property \eqref{goodnormals} holds $\h^2$ a.e. on the specified subset of $\pa \Omega$. In light of this, one could relax the regularity assumption on $\pa \Omega$ to be merely Lipschitz.
\end{Remark}

\begin{figure}
	\centering
	\includegraphics[width = 0.7\textwidth]{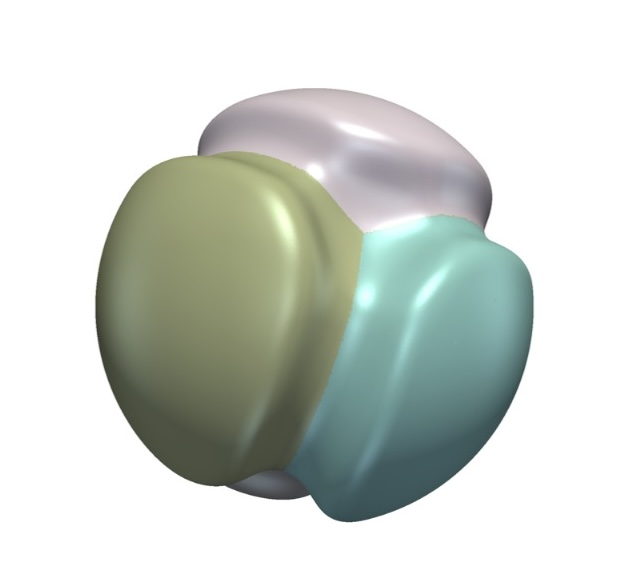}
	\caption{A partition $\ms$ of an admissible domain $\Omega \subset \re^3$.  Here, each of the planes from \cref{Tetrahedral_cone} meets $\pa \Omega$ orthogonally along $\gamma_{ij}$. This property ensures that the constructed $\Omega$ satisfies the condition \eqref{goodnormals}.}
\end{figure}

\section{The main result.}
We now present our main result. In light of previous results on $\Gamma$-convergence this result immediately implies Theorem \ref{pde}. Throughout, we consider a fixed domain $\Omega$ satisfying the conditions \eqref{goodnormals} with respect to the normals coming from $\cone$, and we assume the weights satisfy the condition \eqref{ciclose}.

\begin{Theorem}\label{isolated_local_minimizer}
    The partition $\ms$ is an isolated $L^1$ local minimizer of $E_0$. That is, 
there exists $\delta>0$ such that one has
\begin{equation*}
E_0(\mt) > E_0(\ms)\quad\text{provided}\;0< \lt \left( \mt \Delta \ms\right) \leq \delta.   
\end{equation*}
\end{Theorem}

Before beginning the proof, we consider for any $\delta>0$ the following constrained partition problem:
\begin{equation}
    \inf_{\mt\in\mathcal{A}^\delta}E_0(\mt),\label{conprob}
\end{equation}
where \begin{equation}\label{admissible_class}
    \mathcal{A}^\delta := \{ \mt:=( T_1,T_2,T_3,T_4) \text{ is a partition of }\Omega \text{ s.t } \lt \left( \mt \Delta \ms\right) \leq \delta\}.
    \end{equation}
The existence of a minimizer of this constrained problem is standard via the Direct Method, in light of the lower-semicontinuity of perimeter under $L^1$-convergence. We denote a minimizer of this problem by  
\[\mt^{\delta}= (T_1^\delta,T_2^\delta,T_3^\delta,T_4^\delta).
\] 

The proof of Theorem \ref{isolated_local_minimizer} will rely crucially on the following `corralling' result, which shows that for each $i\in\{1,2,3,4\}$, $\pa T_i\cap\pa\Omega$ lies near  $\pa S_i\cap\pa\Omega$. The argument represents a type of boundary `infiltration lemma', and is an adaptation of the interior result to be found in \cite{Leonardi} and \cite[Lemma 30.2]{Maggi}. 
\begin{Theorem}\label{lemma1}
 There exists a number $\delta>0$ depending on $\cone$ and $\Omega$ such that the following three conditions hold for the partition $\mt^\delta$:
\begin{enumerate}[(I)]
    \item \label{thm4.1}For each $i \in \{ 1,2,3,4\}$, one has :
    \[
    \h^2\bigg(\Big\{ x \in \big( \bigcup_{\substack{j\in\{1,2,3,4\}\\i\not=j}} T_j^\delta\big) \intersect \pa \Omega \intersect \pa S_i : {\rm dist}(x,\cone) \geq \eta \bigg\} \bigg) = 0.
    \]
    \item \label{thm4.2}For distinct indices $i,j,k,l \in \{ 1,2,3,4\}$, one has:
    \[
    \h^2 \left( \{x \;{\color{black}\in \pa \Omega}:\,{\rm dist}(x,q_{ijk})<\eta\} \intersect T_l^\delta\right) = 0.
    \]
    \item\label{thm4.3} For the distinct indices $i,j,k,l \in \{ 1,2,3,4\}$, one has:
    \[
    \h^2 \left( \{x \;{\color{black}\in \pa \Omega }:\,{\rm dist}(x,\gamma_{ij})<\eta\}\intersect \left(T_k^\delta\union T_l^\delta\right)\right) = 0.
    \]
\end{enumerate}

\end{Theorem}

\begin{proof}
Throughout this proof we will suppress the $\delta$ dependence of $\mtd$ and write simply $\mt=(T_1,T_2,T_3,T_4)$. {\color{black} For any $i \in \{1,2,3,4\}$, let $0< r_1<r_0$ be such that for all $x \in \pa S_i \intersect \pa \Omega \text{ with } {\rm dist}(x,\cone)\geq \eta$, one has
\begin{equation*}
B(x,r_1)\intersect \Omega \ssubset S_i,
\end{equation*}
 where the value $r_0$ is defined in \cref{isoperimetric_theorem}.}
We begin with the proof of (I), taking, for instance, the
 case $i=1$; the other cases follow similarly. Consider any $x_0 \in \pa S_1 \intersect\pa \Omega$ such that ${\rm dist}(x_0,\cone) \geq \eta$ and suppose first that for some $0 <r\leq r_1$ one has the condition
    \begin{equation}\label{case1hype}
    \lt \Big(\left[T_2\union T_3 \union T_4\right] \intersect B(x_0,r)\Big) \leq \varepsilon_0 r^3, \text{ where } \varepsilon_0 = \left(\frac{\Lambda}{6}\right)^3,
    \end{equation}
     where $\Lambda = \Lambda (c_1,c_2,c_3,c_4,C_0) >0$ is a constant to be defined later, with $C_0$ defined in \cref{isoperimetric_theorem}.
     Then we will show that 
    \begin{equation}\label{case1conclu}
    \lt\left(\left[T_2\union T_3 \union T_4\right]\intersect B\left(x_0,\frac{r}{2}\right)\right) = 0.
    \end{equation}
To this end, suppressing dependence on $\delta$, we let
\begin{equation}\label{defnms}
m(s) := \lt \big( \left[T_2\union T_3 \union T_4\right] \intersect B(x_0,s)\big).
\end{equation}
Then, in light of the co-area formula, we have
\[
m(s) = \int_{0}^s \sum_{j=2}^4 \h^2 \left( T_j \intersect \pa B(x_0,t)\right)dt.
\]
Further, we have: 
\begin{equation}\label{AbhisDog}
    \begin{aligned}
        m'(s) &= \sum_{j=2}^4 \h^2 \left( T_j \intersect \pa B(x_0,s)\right) \text{ for a.e. } 0<s<r_1,
    \end{aligned}
\end{equation}
and since the perimeter of $T_j$ in $\Omega$ is finite for each $j$, we have
\begin{equation}\label{defnm's2}
        \sum_{j=2}^4 \h^2 \left( \pa^* T_j \intersect \pa B(x_0,s)\right)=0\;\text{ for a.e. } 0<s<r_1.
\end{equation}

For each $0<s<r_1$, we now modify the partition $\mt$ by giving $\left[T_2\union T_3 \union T_4\right] \intersect B(x_0,s)$ to $T_1$ to form a new partition $\mf^s := (F_1^s, F_2^s, F_3^s, F_4^s)$ given by
\begin{equation}\label{modpart}
    \begin{aligned}
        F_1^s &:= T_1 \union \left[B(x_0,s) \intersect \left(T_2\union T_3 \union T_4\right)\right] , \\
        F_2^s &:= T_2 \setminus B(x_0,s), \\
        F_3^s &:= T_3 \setminus B(x_0,s), \\
        F_4^s &:= T_4 \setminus B(x_0,s),\\
    \end{aligned}
\end{equation}
where we again have suppressed the $\delta$ dependence of $\mf^s$ {\color{black} (see \cref{partition_of_T} and \cref{partition_of_F})}. We first observe that the partition $\mf^s$ at worst leaves the $L^1$ distance to $\ms$ unchanged; hence the new partition satisfies:
    \[
    \lt (\mf^s \Delta \ms) \leq \delta \quad \text{for all}\; 0<s<r_1.
    \]
Thus, $\mf^s$ is an element of $\mathcal{A}^\delta$. Since $\mt = (T_1,T_2,T_3,T_4)$ is the minimizer of \eqref{conprob}, we have the minimality condition:
\begin{equation}\label{minimality_1}
E_0(\mt)\leq E_0(\mf^s)\quad\text{for all}\;s\in (0,r_1).
\
\end{equation}

\begin{figure}
    \centering
    \includegraphics[width=0.7\linewidth]{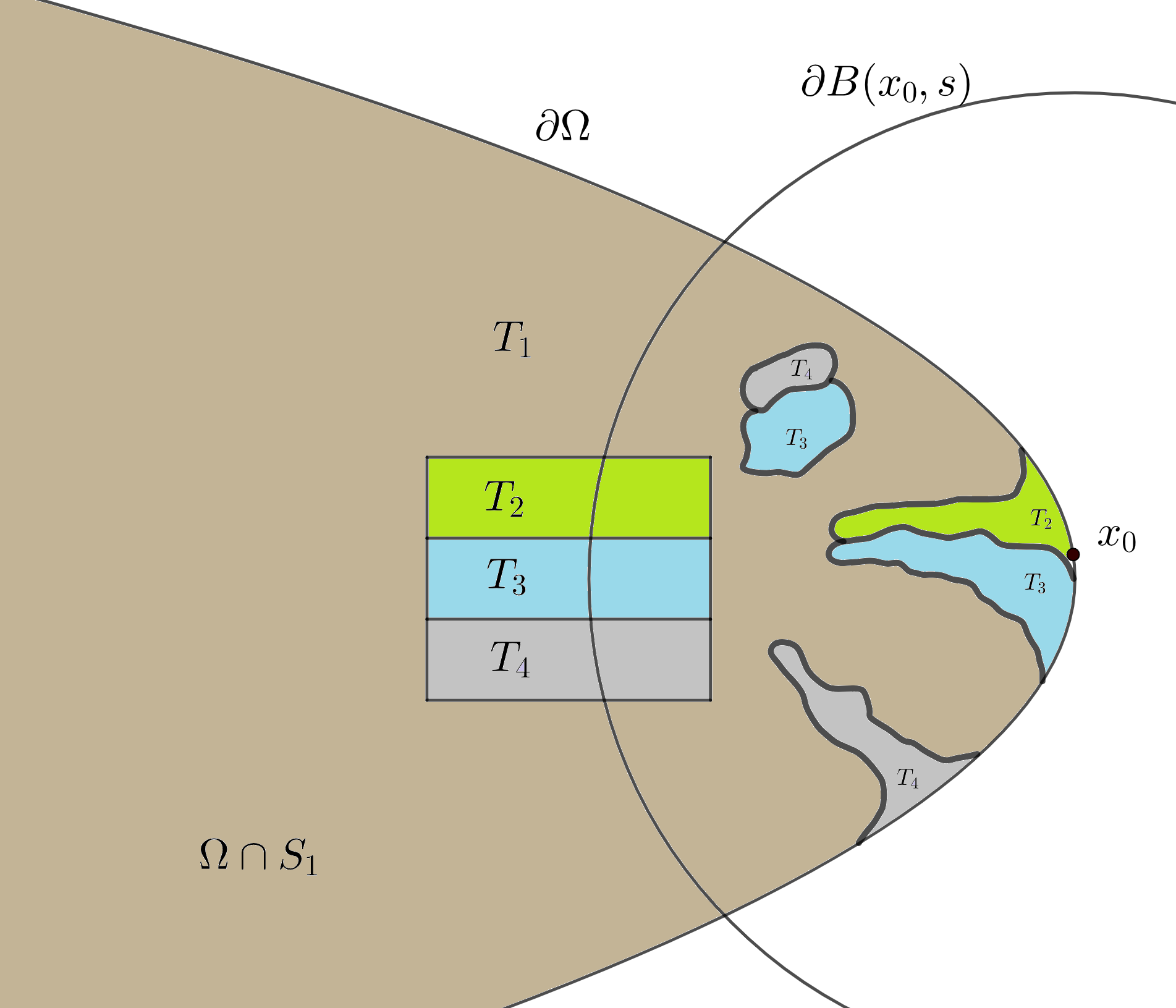}
    \caption{The partition $\mt = (T_1,T_2,T_3,T_4)$ in $\Omega \intersect S_1$.}
    \label{partition_of_T}
\end{figure}

\begin{figure}
    \centering
    \includegraphics[width=0.7\linewidth]{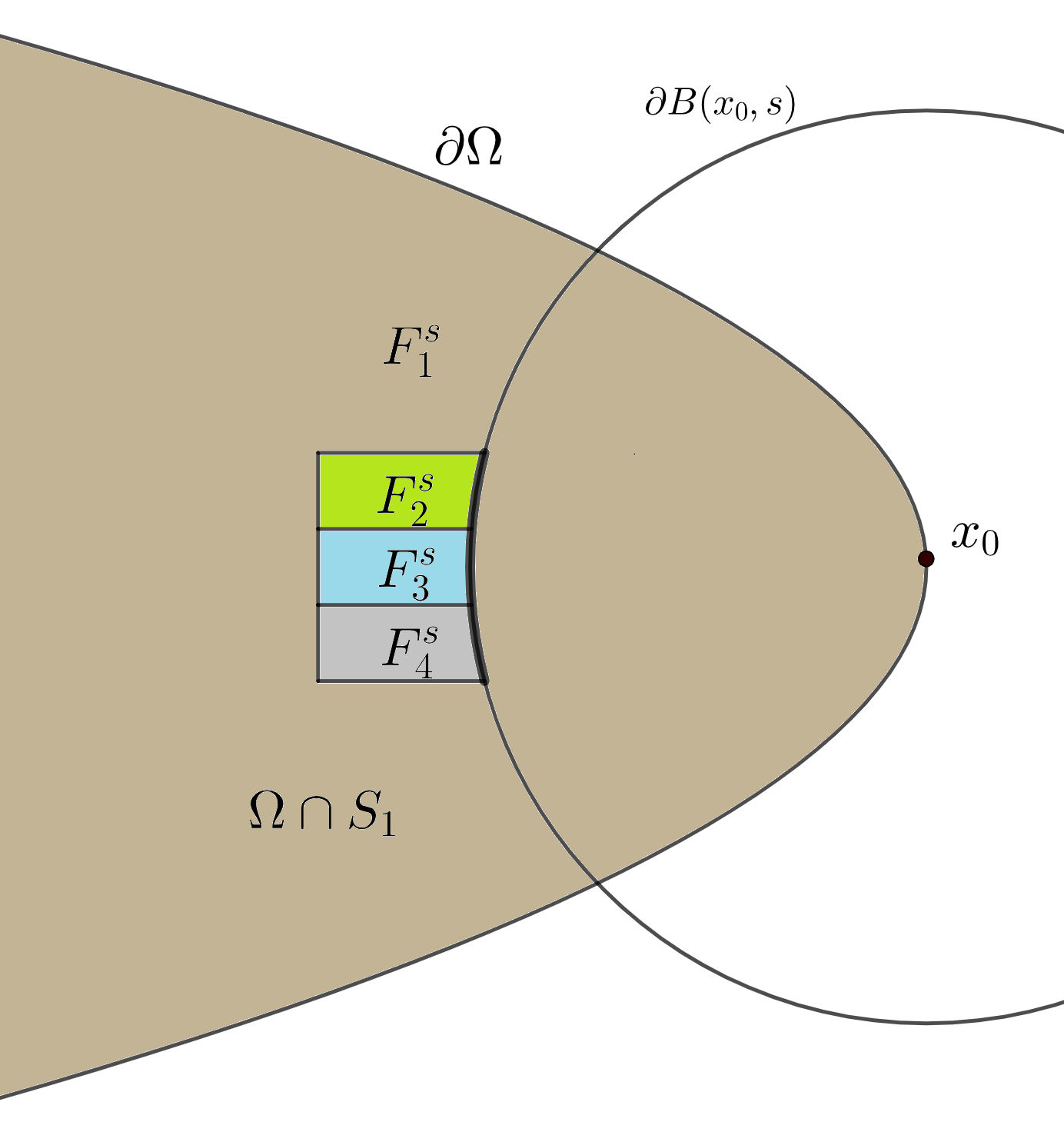}
    \caption{The partition $\mf^s = (F_1^s,F_2^s,F_3^s,F_4^s)$ in $\Omega  \intersect S_1$.}
    \label{partition_of_F}
\end{figure}

At this point, we introduce the simplifying notation
\begin{equation}\label{Tij}
\begin{aligned}
    T_{ij}^s&:=\h^2\big(\pa^*T_i\cap\pa^*T_j\cap B(x_0,s)\big).
    \end{aligned}
\end{equation}
We also set 
\[
c_{max} := \max_{1 \leq j \leq 4} c_j \text{ and }c_{min} := \min_{1 \leq j \leq 4} c_j.
\]
The main idea of this step is to express the terms on either side of the inequality in \eqref{minimality_1} using \eqref{ciE} to get a differential inequality. To do so, we first express $E_0(\mt)$ using the notation \eqref{ciE} and the equation \eqref{defnm's2} as:
 \begin{equation}\label{expand_mtd}
     \begin{aligned}
          E_0(\mt)
          &=  \sum_{j=1}^4 c_j \h^2 (\pa^* T_j \intersect \Omega)\\
          &=\sum_{j=1}^4 c_j \h^2 (\pa^* T_j \intersect \Omega \intersect B(x_0,s)) + \sum_{j=1}^4 c_j \h^2 (\pa^* T_j \intersect \Omega \intersect \pa B(x_0,s)) \\
          &\qquad \qquad + \sum_{j=1}^4 c_j \h^2 (\pa^* T_j \intersect \Omega \intersect [\Omega \setminus \overline{B(x_0,s)}])\\
          & = \Big[[c_1 + c_2]T_{12}^s + [c_1 + c_3]T_{13}^s + [c_1 + c_4]T_{14}^s + [c_2 + c_3]T_{23}^s \\
          &\qquad  + [c_2 + c_4]T_{24}^s + [c_3 + c_4]T_{34}^s\Big] + \sum_{j=1}^4 c_j \h^2 (\pa^* T_j \intersect [\Omega \setminus \overline{B(x_0,s)}]) \\
          &\geq 2c_{min} \left[ T_{12}^s + T_{13}^s + T_{14}^s\right ]+ \sum_{j=1}^4 c_j \h^2 (\pa^* T_j    \intersect [\Omega \setminus \overline{B(x_0,s)}])  \\
          & = 2c_{min} {\rm Per}\big[ ( T_2 \union T_3 \union T_4);\Omega \intersect B(x_0,s) \big] + \sum_{j=1}^4 c_j \h^2 (\pa^* T_j \intersect [\Omega \setminus \overline{B(x_0,s)}]), 
    \end{aligned}
 \end{equation}
Then in light of \eqref{defnm's2} and \eqref{modpart}, one observes that for all $0<s<r_1$,
\begin{equation*}
    \sum_{j=1}^4
{\rm Per}(F_j^s;\Omega \intersect B(x_0,s)) = 0.
\end{equation*}
Hence, recalling \eqref{AbhisDog}, we obtain
\begin{equation}\label{expand_mfd}
\begin{aligned}
          E_0(\mf^s) = & \sum_{j=1}^4 c_j \h^2 (\pa^* F_j^s \intersect \Omega)\\
          & = [c_1 +c_2]\h^2( T_2 \intersect\pa B(x_0,s)) + [c_1 + c_3] \h^2(T_3 \intersect \pa B(x_0,s)) \\
          &+ [c_1 + c_4] \h^2 ( T_4 \intersect \pa B(x_0,s)) + \sum_{j=1}^4 c_j \h^2 (\pa^* F_j^s\intersect [\Omega \setminus \overline{B(x_0,s)}])\\
          &\leq 2c_{max} m'(s) + \sum_{j=1}^4 c_j \h^2 (\pa^* F_j^s\intersect [\Omega \setminus \overline{B(x_0,s)}]).
\end{aligned}
\end{equation}

We observe that outside the ball $B(x_0,r)$, the partitions $\mt$ and $\mf^s$ are the same. Thus, substituting \eqref{expand_mtd} and \eqref{expand_mfd} in the minimality condition \eqref{minimality_1}, one has 
\begin{equation*}
 c_{min} {\rm Per} \big[(T_2 \union T_3 \union T_4) ; \Omega \intersect B(x_0,s) \big] \leq c_{max} m'(s). 
\end{equation*}
Now we add $c_{min}\sum_{j=2}^4 \h^2\big(T_j \intersect \pa B(x_0,s)\big)=c_{min}m'(s)$ to both sides of this inequality to obtain
\begin{equation*}
        c_{min}{\rm Per} \big[(T_2 \union T_3 \union T_4)\intersect B(x_0,s) ; \Omega  \big] \leq (c_{min} + c_{max}) m'(s).        
\end{equation*}
Applying the isoperimetric inequality (\cref{isoperimetric_theorem}) to the left-hand side, we arrive at the differential inequality

\begin{equation}\label{lambda_1}
    \Lambda m(s)^{\frac{2}{3}}\leq m'(s), 
\end{equation}
where $\Lambda$ is any constant satisfying
\begin{equation}\label{Prathibha}
\Lambda \leq \left(\frac{C_0 \; c_{min}}{c_{max} + c_{min}}\right).
\end{equation}

Now, since the map $s \mapsto m(s)$ is monotonically increasing, we have that the support of $m$ is $[r^*,\infty)$, for some $r^* \geq 0$. If $m(r) = 0$ for $r$ in \eqref{case1hype}, then trivially \eqref{case1conclu} holds. Otherwise $m(r)>0$ and so necessarily $r > r^*$. Integrating \eqref{lambda_1} on $(r^* ,r)$ gives
\begin{equation}\label{rr*}
\begin{aligned}
    r - r^* \leq \frac{3}{\Lambda} \left[ \sqrt[3]{m(r)} - \sqrt[3]{m(r^*)}\right]& = \frac{3}{\Lambda}\sqrt[3]{m(r)} \\
    &\leq\frac{3}{\Lambda} \sqrt[3]{\varepsilon_0}r \leq \frac{r}{2},
    \end{aligned}
\end{equation}
where we used \eqref{case1hype} for the last two inequalities in \eqref{rr*}. Hence,  $r^*\geq r/2>0$, and so the assumption \eqref{case1hype} implies the conclusion \eqref{case1conclu}.

On the other hand, consider any $x_0 \in \pa S_1 \intersect \pa \Omega$ such that ${\rm dist}(x_0,\cone) \geq \eta$ for which assumption \eqref{case1hype} does not hold. That is, suppose
    \begin{equation}\label{counterlem1}
      \forall r \leq r_1, \quad  \lt(\left[T_2\union T_3 \union T_4\right]\intersect B(x_0,r)) > \varepsilon_0 r^3.
    \end{equation}
 In particular, take $r = \sqrt[4]{\delta}$, where we stipulate that 
 \begin{equation}\label{Karthik}
 \delta < \min\left\{\left(\frac{r_1}{2}\right)^4, \frac{\varepsilon_0^4}{2}\right\}.
 \end{equation}
 Plugging this $r$ in \eqref{counterlem1} and noting that $\lt(\mt \Delta \ms) \leq \delta$, we reach a contradiction to \eqref{Karthik} in that 
\begin{equation*}
    \delta \geq \lt(\left[T_2\union T_3 \union T_4\right] \intersect B(x_0,r_1)) > \varepsilon_0 r^3 = \varepsilon_0 \delta^{\frac{3}{4}},
\end{equation*}
and so \eqref{counterlem1} cannot hold. The proof of (I) is complete.

Now, we turn to the proof of (II). We will consider the case $i=1,j=2,k=3,l=4$, so our goal is to show
\begin{equation*}
  \h^2\left (T_4 \intersect \{ x \in \pa \Omega : \; {\rm dist}(x,q_{123})<\eta\} \right) = 0.  
\end{equation*}
The other cases follow in the similar manner. To this end, we fix any $0 < r_2 < r_0$ such that
\begin{equation*}
r_2 < \frac{1}{2} \min\left\{ {\rm dist}\left( B_{i_1i_2i_3}, B_{j_1j_2j_3}\right) \left| \begin{aligned}
    &i_1,j_1,i_2,j_2,i_3,j_3 \in \{1,2,3,4\},\\
    & i_k \neq j_k \text{ for all }k=1,2,3.
    \end{aligned}\right.\right\}
\end{equation*}
where $ B_{ijk} := \{x \in \pa \Omega ;\; {\rm dist}(x,q_{ijk}) < \eta\}$ and again, 
the value $r_0$ comes from \cref{isoperimetric_theorem}. 

We consider a point $x_0 \in B_{123}$. In analogy with \eqref{case1hype} from the proof of (I), we assume that for some $0<r<r_2$, one has
\begin{equation}\label{X1}
\lt \left( B(x_0,r) \intersect T_4\right) \leq \varepsilon_0 r^3, \text{ with } \varepsilon_0=\left(\frac{\Lambda}{6}\right)^3,
\end{equation}
where $\Lambda$ is to be specified later. A crucial step here is that by an application of the pigeon hole principle, we see that for every $0<s<r$, there exists $j_0 = j_0(s) \in \{1,2,3\}$ such that one has
\begin{equation*}
T_{j_0 4}^s \geq\frac{1}{3} \left[ T_{14}^s + T_{24}^s + T_{34}^s\right],
\end{equation*}
or equivalently,
\begin{equation}\label{pigeon_hole}
    2T_{j_04}^s\geq \sum_{\substack{j \in \{1,2,3\} \\ j \neq j_0}} T_{j4}^s
\end{equation}
in the notation of \eqref{Tij}.
The primary difference between the proof of (II) and (I) is that here we create a competing partition by giving $ B(x_0,s)\intersect T_4$ to $T_{j_0}$, rather than giving all of $B(x_0,s) \intersect (\union_{j=2}^4 T_j)$ to $T_{j_0}$. More precisely, for all $0 < s < r$, we set
    \begin{equation*}
        \begin{aligned}
            G_{j_0}^s &:= T_{j_0}\union \Big[ B(x_0,s)\intersect T_4 \Big], \\
            G_k^s &:= T_k, \text{ if } k \in\{1,2,3\}, k \neq j_0 \\
            G_4^s &:= T_4 \setminus  B(x_0,s), \\
            \mathcal{G}^s &:= (G_1^s,G_2^s,G_3^s,G_4^s),
        \end{aligned}
    \end{equation*}
    where we again have suppressed the $\delta$ dependence of $\mathcal{G}^s$. Since in this case, we are working in a neighborhood of $q_{123}$, giving $T_4 \intersect B(x_0,s)$ to $T_j$ either decreases the $L^1$ distance to $\ms$ or at worst, the $L^1$ distance remains unchanged. Therefore, we have
\[ 
 \lt \Big( \ms \Delta \mathcal{G}^s\Big)\leq \delta \text{ for all } 0<s<r \text{ and so } \mathcal{G}^s \in \mathcal{A}^\delta,\]
thus implying the minimality condition:
\begin{equation}\label{minimality_case_2}
    E_0(\mt) \leq E_0(\mathcal{G}^s).
\end{equation}
Here, the role of $m(s)$ from \eqref{defnms} is taken up by say $n(s)$, where for all $0<s<r$, we set
    \begin{equation*}
        \begin{aligned}
    n(s) &:= \lt (T_4 \intersect  B(x_0,s)) \\
    &= \int_0^s \h^2 \Big[ T_4 \intersect \pa B(x_0,t)\Big] dt.
        \end{aligned}
    \end{equation*}
    Further, for a.e. $0<s<r$, we have:
    \begin{equation}\label{Adimurthi}
        \begin{aligned}
           n'(s) &= \h^2 \Big[T_4 \intersect \pa B(x_0,s) \Big] \text{ and } \h^2 \Big[\pa^* T_4 \intersect \pa B(x_0,s)\Big]=0,
        \end{aligned}
    \end{equation}
where again we have suppressed the $\delta$ dependence of $n$. 

Let us consider an $s$-value in $(0,r_2)$ where $j_0 =1$. To this end, \eqref{minimality_case_2} takes the form
 \begin{equation*}
     \begin{aligned}
          &E_0(\mt) =  \sum_{j=1}^4 c_j \h^2 (\pa^* T_j \intersect \Omega)\\
          &= \Big[[c_1 + c_2]T_{12}^s + [c_1 + c_3]T_{13}^s + [c_1 + c_4]T_{14}^s + [c_2 + c_3]T_{23}^s \\
          &\qquad  + [c_2 + c_4]T_{24}^s + [c_3 + c_4]T_{34}^s\Big] + \sum_{j=1}^4 c_j \h^2 (\pa^* T_j \intersect [\Omega \setminus \overline{B(x_0,s)}])\\
          &\leq  E_0(\mathcal{G}^s) = \sum_{j=1}^4 c_j \h^2 (\pa^* G_j^s \intersect \Omega)\\
          & = [c_1+c_2](T_{12}^s +T_{24}^s) + [c_1 + c_3](T_{13}^s + T_{34}^s) +[c_2 + c_3] T_{23}^s \\
          &+ [c_1 + c_4] \h^2(T_4 \intersect \pa B(x_0,r)) +\sum_{j=1}^4 c_j \h^2 (\pa^* F_j^s\intersect [\Omega \setminus \overline{B(x_0,r)}]).
\end{aligned}
\end{equation*}

Consequently, after cancellation of like terms, we have
\begin{equation*}
    [c_1 + c_4] T_{14}^s  + c_4(T_{24}^s + T_{34}^s) \leq c_1 (T_{24}^s + T_{34}^s) + [c_1 + c_4]n'(s),
\end{equation*}
where we recall \eqref{Adimurthi}. Adding $c_4 n'(s)$ to both the sides and noting that ${ \rm Per}(T_4 \intersect B(x_0,s); \Omega ) = T_{14}^s + T_{24}^s + T_{34}^s + n'(s)$, we see that
\begin{equation}\label{Jura}
    c_1 T_{14}^s + c_4 {\rm Per}(T_4 \intersect B(x_0,s); \Omega) \leq c_1 (T_{24}^s + T_{34}^s) + [c_1 + 2c_4] n'(s).
\end{equation}

Now, we split the argument into two cases. In light of \eqref{pigeon_hole}, either one has
\begin{equation}\label{Adi1}
    T_{14}^s \geq \frac{1}{2} \left(T_{14}^s + T_{24}^s + T_{34}^s\right)
\end{equation}
or else 
\begin{equation}\label{Adi2}
    \frac{1}{3} \left(T_{14}^s + T_{24}^s + T_{34}^s\right) \leq T_{14}^s < \frac{1}{2} \left(T_{14}^s + T_{24}^s + T_{34}^s\right).
\end{equation}

Pursuing the case of \eqref{Adi1}, we have
\[
T_{14}^s \geq T_{24}^s + T_{34}^s,
\]
which we apply to \eqref{Jura} to obtain
\[
c_4 {\rm Per}(T_4 \intersect B(x_0,s); \Omega) \leq  [c_1 + 2c_4] n'(s).
\]
Utilizing the isoperimetric inequality, we arrive at a differential inequality
\[
n'(s) \geq \left( \frac{c_4 \,C_0 }{ c_1 + 2c_4}\right) n(s)^{\frac{2}{3}}.
\]

More generally, for any value of $j_0\in\{1,2,3\}$ we arrive at the inequality
\begin{equation*}
    \begin{aligned}
        n'(s) \geq \left( \frac{c_4 \,C_0 }{ c_{j_0} + 2c_4}\right) n(s)^{\frac{2}{3}} &\geq \left( \frac{c_4 \,C_0 }{ c_{max} + 2c_4}\right) n(s)^{\frac{2}{3}}\\
        &\geq \left( \frac{c_{min} \,C_0 }{ 3 c_{max} }\right) n(s)^{\frac{2}{3}}
    \end{aligned}
\end{equation*}

We now pursue the case \eqref{Adi2} for $j_0 = 1$. Using $ 2T_{14}^s \geq T_{24}^s + T_{34}^s>T_{14}^s$ in \eqref{Jura}, we now find that 
 \begin{equation*}
 \begin{aligned}
  &c_4 {\rm Per}(T_4 \intersect B(x_0,s); \Omega)\\ 
  &\leq c_1 T_{14} + (c_1 + 2c_4) n'(s)\\
  &\leq \frac{c_1}{2} \Big(T_{14}^s + T_{24}^s + T_{34}^s + \h^2(T_4\intersect\pa B(x_0,s)\Big) + \Big[ \frac{c_1}{2} + 2c_4\Big] n'(s)\\
  & = \frac{c_1}{2} {\rm Per}(T_4 \intersect B(x_0,s);\Omega) + \Big[ \frac{c_1}{2} + 2c_4\Big] n'(s).
  \end{aligned}
 \end{equation*}
Combining like terms  and applying the isoperimetric inequality, we obtain the differential inequality
\[
n'(s) \geq \left( \frac{[2c_4 -c_1] C_0}{c_1 + 4c_4}\right) n(s)^{\frac{2}{3}},
\]
where the coefficient on the right is strictly positive in light of our assumption \eqref{ciclose}.
More generally, for any value of $j_0(s)\in\{1,2,3\}$ we arrive at the inequality
\begin{equation*}
    \begin{aligned}
        n'(s)  &\geq \left( \frac{[2c_4 -c_{j_0}] C_0}{c_{j_0} + 4c_4}\right) n(s)^{\frac{2}{3}} \\
        & \geq \left( \frac{\min_{i\neq j}(2c_i - c_j) C_0}{5c_{max}}\right)n(s)^{\frac{2}{3}},
    \end{aligned}
\end{equation*}
where the hypothesis \eqref{ciclose} guarantees that the coefficient on the right is strictly positive. Therefore, for any 
\begin{equation}\label{Prathibha2}
\Lambda < \min \left\{ \left( \frac{c_{min} \,C_0 }{ 3 c_{max} }\right),\left( \frac{\min_{i\neq j}(2c_i - c_j) C_0}{5c_{max}}\right)\right\}, 
\end{equation}
we obtain the differential inequality
\[
\Lambda n(s)^{\frac{2}{3}} \leq n'(s).
\]
Clearly, the same differential inequality will hold for any $\eta-$neighborhood of $q_{ijk}$. The rest of the argument follows exactly as in case (I). 

When the reverse inequality in \eqref{X1} holds, we can choose 
\[
\delta < \min\left\{\left(\frac{r_2}{2}\right)^4, \frac{\varepsilon_0^4}{2}\right\}
\]
to conclude the proof for case (II).

We now turn to the proof of (III) and consider, for example, the case $i=1,j=2,k=3,l=4$; the other cases follow similarly. That is, we need to show that 
\begin{equation*}
   \h^2 \left( \{x \in \pa \Omega :\; {\rm dist}(x,\gamma_{12})<\eta\} \intersect (T_3\union T_4)\right) = 0. 
\end{equation*}
Proceeding as in the proofs of (I) and (II), we let $0 < r_3<r_0$ be a number such that 
\begin{equation*}
\{x \in \pa \Omega :\; {\rm dist}(x,\gamma_{12})<\eta\} \intersect B(x,r_3) \ssubset \left(\pa S_1 \union \pa S_2\right)\intersect \pa \Omega.
\end{equation*}

 Now, let $x_0 \in \{x \in \pa \Omega :\; {\rm dist}(x,\gamma_{12})<\eta\}$ and assume that for some $0<r<r_3$, one has
\begin{equation}\label{X}
\lt((T_3\union T_4) \intersect B(x_0,r)) \leq \varepsilon_0 r^3, \text{ with } \varepsilon_0 := \left(\frac{\Lambda}{6}\right)^3,
\end{equation}
where the $\Lambda$ from \eqref{Prathibha2} will suffice.
We first observe that if we consider a new partition in which, say, $T_4\cap B(x_0,s)$ is given any $T_j$ for $j\in\{1,2,3\}$, then at worst, then the constraint condition in \eqref{admissible_class} is unchanged. This puts us in exactly the setting of case (II). Hence, by the same analysis, we may immediately conclude that 
$T_4\cap B(x_0,r/2)$ is empty. 

It remains to argue that $T_3\cap B(x_0,s)$ is empty for small $s$ as well. This scenario is precisely the one considered in \cite{MikeLia}, but for the sake of completeness, we will provide the proof. 

To this end, there are two possibilities: either one has $T_{13}^s \geq T_{12}^s$ or that $T_{12}^s > T_{13}^s$. We pursue the first possibility; the arguments are the same for the other case. We begin by giving $T_3 \intersect B(x_0,s)$ to $T_1$ and we call this new partition $\h^s$. We note that $\h^s \in \mathcal{A}^\delta$.

Using the minimality condition $E_0(\mt) \leq E_0(\h^s)$ we see that
\begin{align}
   & c_1(T_{12}^s +T_{13}^s) + c_2 (T_{12}^s +T_{23}^s) + c_3 (T_{13}^s + T_{23}^s)\nonumber \\
   &\leq [c_1 + c_2] (T_{12}^s + T_{23}^s) +[c_1 + c_3] \h^2 (T_3 \intersect \pa B(x_0,s)).\nonumber
\end{align}

Canceling like terms, using $T_{13}^s \geq T_{23}^s$ and adding $c_3  \h^2 (T_3 \intersect \pa B(x_0,s))$ to both sides of the resulting inequality, one obtains
\[
c_3 {\rm Per}(T_3\intersect B(x_0,s);\Omega ) \leq [c_1 + 2c_3]  \h^2 (T_3 \intersect \pa B(x_0,s)).
\]

Introducing $p(s) := \lt( T_3 \intersect B(x_0,s))$, applying the isoperimetric inequality, one arrives at the differential inequality
\[
\left(\frac{c_3 C_0}{c_1 + 2 c_3}\right) p(s)^{\frac{2}{3}} \leq p'(s).
\]

If $T_{12}^s > T_{13}^s$, then by giving $T_3 \intersect B(x_0,s)$ to $T_2$ instead of $T_1$, the only change in the inequality is that $c_1$ is replaced by $c_2$ in the denominator. Arguing as before, working with the differential inequality 
\begin{equation}\label{Prathibha3}
\Lambda p(s)^{\frac{2}{3}} \leq p'(s) \text{ for any } \Lambda \leq \frac{c_{min}C_0}{3c_{max}},
\end{equation}
we conclude that $p(r/2) = 0$ and so there is no $T_3$ in $B(x_0,r/2)$ is empty under the assumed inequality \eqref{X}. For the reverse inequality, we arrive at the same conclusion by choosing
\[
\delta < \min\left\{\left(\frac{r_3}{2}\right)^4, \frac{\varepsilon_0^4}{2}\right\}.
\]

Finally, for $\Lambda$ satisfying \eqref{Prathibha}, \eqref{Prathibha2} and \eqref{Prathibha3}, \cref{lemma1} is established provided
\[
\delta < \min\left\{\left(\frac{r_1}{2}\right)^4, \left(\frac{r_2}{2}\right)^4, \left(\frac{r_3}{2}\right)^4, \frac{1}{2}\left(\frac{\Lambda}{6}\right)^4\right\}.
\]

\end{proof}

Having argued that the constrained minimizer $\mtd$ to problem \eqref{conprob} has the property that for each $i\in\{1,2,3,4\}$, $\pa T_i$ lies in an $\eta$ neighborhood of $\pa S_i$  on $\pa\Omega$, we can now apply a calibration argument inspired by that of \cite{Morgan}.

\begin{proof}[Proof of \cref{isolated_local_minimizer}] 
With \cref{lemma1} in hand, we now argue that in fact, the solution $\mtd$ to \eqref{conprob} coincides with $\ms$. In what follows, for ease of notation, we will again not indicate the $\delta$ dependence of $\mtd$ and write simply $\mt=(T_1,T_2,T_3,T_4)$. We recall that $n_{ij}$ denotes the unit normal to $\pa S_i \intersect \pa S_j$ pointing from $S_i$ to $S_j$, $\nu_{\pa \Omega}$ denotes the unit outward normal to $\pa\Omega$, and we let $\nu_{ij}^*$ denote the unit normal to $\pa^* T_i \intersect \pa^* T_j$ pointing from $T_i$ to $T_j$, 

Throughout this section, we use the formulation \eqref{originalE0} of $E_0$ rather than \eqref{ciE}. We begin with four applications of the Gauss-Green theorem to obtain
\begin{equation}\label{Q23}
    \begin{aligned}
        0 &= \int_{T_1} {\rm div}(n_{14}) dx = \int_{\pa^*T_1 \intersect \pa^* T_2}n_{14}\cdot \nu_{12}^* \,d\h^2 + \int_{\pa^*T_1 \intersect \pa^* T_3}n_{14}\cdot \nu_{13}^* \,d\h^2 \\
        & \quad \quad \quad \quad \quad \quad \quad \quad +\int_{\pa^*T_1 \intersect \pa^* T_4}n_{14}\cdot \nu_{14}^* \,d\h^2 +  \int_{ \pa^*T_1 {\color{black}\intersect \pa \Omega}} n_{14}\cdot \nu_{\pa\Omega} \,d\h^2, 
    \end{aligned}
\end{equation}

\begin{equation}\label{Q22}
    \begin{aligned}
        0 &= \int_{T_2} {\rm div}(n_{14}) dx = \int_{\pa^*T_2 \intersect \pa^* T_1}n_{14}\cdot \nu_{21}^* \,d\h^2 + \int_{\pa^*T_2 \intersect \pa^* T_3}n_{14}\cdot \nu_{23}^* \,d\h^2 \\
        & \quad \quad \quad \quad \quad \quad \quad \quad +\int_{\pa^*T_2 \intersect \pa^* T_4}n_{14}\cdot \nu_{24}^* \,d\h^2 +  \int_{ \pa^*T_2 {\color{black}\intersect \pa \Omega}}n_{14}\cdot \nu_{\pa \Omega} \,d\h^2,
    \end{aligned}
\end{equation}
\begin{equation}\label{Q21}
    \begin{aligned}
        0 &= \int_{T_2} {\rm div}(n_{21}) dx = \int_{\pa^*T_2 \intersect \pa^* T_1}n_{21}\cdot \nu_{21}^* \,d\h^2 + \int_{\pa^*T_2 \intersect \pa^* T_3}n_{21}\cdot \nu_{23}^* \,d\h^2 \\
        & \quad \quad \quad \quad \quad \quad \quad \quad +\int_{\pa^*T_2 \intersect \pa^* T_4}n_{21}\cdot \nu_{24}^* \,d\h^2 +  \int_{ \pa^*T_2 {\color{black}\intersect \pa \Omega}}n_{21}\cdot \nu_{\pa \Omega} \,d\h^2,
    \end{aligned}
\end{equation}

and

\begin{equation}\label{Q24}
    \begin{aligned}
        0 &= \int_{T_3} {\rm div}(n_{34}) dx = \int_{\pa^*T_3 \intersect \pa^* T_1}n_{34}\cdot \nu_{31}^* \,d\h^2 + \int_{\pa^*T_3 \intersect \pa^* T_2}n_{34}\cdot \nu_{32}^* \,d\h^2 \\
        & \quad \quad \quad \quad \quad \quad \quad \quad +\int_{\pa^*T_3 \intersect \pa^* T_4}n_{34}\cdot \nu_{34}^* \,d\h^2 +  \int_{ \pa^*T_3 {\color{black}\intersect \pa \Omega}} n_{34}\cdot \nu_{\pa \Omega} \,d\h^2.
    \end{aligned}
\end{equation}
We remark that the choice of the normals and the domains of integration is not unique. For example, we could also have also have begun with the four integrals
\[ \int_{T_4} {\rm div}(n_{41})\, dx,\;  \int_{T_4}{\rm div} (n_{13})\,dx ,\; \int_{T_1} {\rm div}(n_{13})\,dx \;\mbox{and}\; \int_{T_2}{\rm div}(n_{23})\,dx.
\]

Multiplying ${\rm div} (n_{ij})$ by $c_{ij}$ in  equations \eqref{Q23}, \eqref{Q22}, \eqref{Q21} and \eqref{Q24}, we add them and then use the cyclic property of the normals \eqref{gencij}, the fact that $\nu_{ij}^* = - \nu_{ji}^*$ and $n_{ij} = - n_{ji}$ to obtain 
\begin{equation*}
    \begin{aligned}
        0 &= \int_{T_1} c_{14}{\rm div}(n_{14})\,dx + \int_{T_2}c_{14} {\rm div}(n_{14}) dx \\
        & \qquad \qquad \qquad  + \int_{T_2} c_{12}{\rm div}(n_{21})\,dx + \int_{T_3} c_{34}{\rm div}(n_{34})\,dx\\
        &= \sum_{1 \leq i < j \leq 4} \int_{\pa^*T_i \intersect \pa^*T_j} c_{ij} n_{ij} \cdot \nu_{ij}^* \, d\h^2 - \sum_{i\in\{1,2,3\}}\int_{ \pa^*T_i {\color{black}\intersect \pa \Omega}}c_{i4}n_{4i}\cdot \nu_{\pa\Omega} \,d\h^2.
    \end{aligned}
\end{equation*}
Rearranging this identity and noting that  $n_{ij}\cdot \nu_{ij}^* \leq 1$ for all $i,j$, we have
\begin{equation}\label{1_final}
    \begin{aligned}
        E_0(\mt)&:= \sum_{1\leq i <j \leq 4} c_{ij}\h^2 \left( \pa^* T_i \intersect \pa^* T_j \intersect \Omega\right) \\
         &\geq \sum_{1 \leq i < j \leq 4} \int_{\pa^*T_i \intersect \pa^*T_j} c_{ij} n_{ij} \cdot \nu_{ij}^* \, d\h^2 = \sum_{i\in\{1,2,3\}}\int_{ \pa^*T_i {\color{black}\intersect \pa \Omega}}c_{i4}n_{4i}\cdot \nu_{\pa\Omega} \,d\h^2.
    \end{aligned}
\end{equation}
Here we also observe that if we replace $\mt$ by $\ms$ in the calculation above, then one instead obtains the identity 
\begin{equation}\label{E_0_s}
E_0(\ms) = \sum_{i\in\{1,2,3\}}\int_{ \pa^*S_i {\color{black}\intersect \pa \Omega}}c_{i4}n_{4i}\cdot \nu_{\pa\Omega} \,d\h^2. 
\end{equation}
We will first argue that $\mt$ and $\ms$ agree on $\pa\Omega$.  To this end, we decompose the domains of integration in the last line of \eqref{1_final}, by writing each $\pa\Omega\cap\pa^*T_i$ as $\pa\Omega\cap\pa S_i$ plus corrections, for $i=1,2,3$:

\begin{equation*}
    \begin{aligned}
        &E_0(\mt) \geq \sum_{i\in\{1,2,3\}}\int_{ \pa^*T_i {\color{black}\intersect \pa \Omega}}c_{i4}n_{4i}\cdot \nu_{\pa\Omega} \,d\h^2 = \sum_{i\in\{1,2,3\}}\int_{\pa^*S_i {\color{black}\intersect \pa \Omega}}c_{i4}n_{4i}\cdot \nu_{\pa\Omega} \,d\h^2 \\
        & + \int_{(\pa^* T_2 \union \pa^* T_3 \union \pa^* T_4) \intersect \pa S_1 {\color{black}\intersect \pa \Omega}}c_{14} n_{14}\cdot \nu_{\pa \Omega} \,d\h^2 + \int_{\pa^* T_1 \intersect (\pa S_2 \union \pa S_3 \union \pa S_4) {\color{black}\intersect \pa \Omega}} c_{14}n_{41}\cdot \nu_{\pa \Omega} d \h^2 \\
         &  +\int_{(\pa^* T_1 \union \pa^* T_3 \union \pa^* T_4) \intersect \pa S_2 {\color{black}\intersect \pa \Omega}} c_{24}n_{24}\cdot \nu_{\pa \Omega} \,d\h^2 +\int_{\pa^* T_2 \intersect (\pa S_1 \union \pa S_3 \union \pa S_4) {\color{black}\intersect \pa \Omega}} c_{24}n_{42}\cdot \nu_{\pa \Omega} d \h^2\\
        &  + \int_{(\pa^* T_1 \union \pa^* T_2 \union \pa^* T_4) \intersect \pa S_3 {\color{black}\intersect \pa \Omega}} c_{34}n_{34}\cdot \nu_{\pa \Omega} \,d\h^2 + \int_{\pa^* T_3 \intersect (\pa S_1 \union \pa S_2 \union \pa S_4) {\color{black}\intersect \pa \Omega}} c_{34}n_{43}\cdot \nu_{\pa \Omega} d \h^2.
    \end{aligned}
\end{equation*}

By \eqref{E_0_s}, we have
\begin{equation*}
    \begin{aligned}
        &E_0(\mt) \geq E_0(\ms) \\
        & + \int_{(\pa^* T_2 \union \pa^* T_3 \union \pa^* T_4) \intersect \pa S_1 {\color{black}\intersect \pa \Omega}}c_{14} n_{14}\cdot \nu_{\pa \Omega} \,d\h^2 + \int_{\pa^* T_1 \intersect (\pa S_2 \union \pa S_3 \union \pa S_4) {\color{black}\intersect \pa \Omega}} c_{14}n_{41}\cdot \nu_{\pa \Omega} d \h^2 \\
         &  +\int_{(\pa^* T_1 \union \pa^* T_3 \union \pa^* T_4) \intersect \pa S_2 {\color{black}\intersect \pa \Omega}} c_{24}n_{24}\cdot \nu_{\pa \Omega} \,d\h^2 +\int_{\pa^* T_2 \intersect (\pa S_1 \union \pa S_3 \union \pa S_4) {\color{black}\intersect \pa \Omega}} c_{24}n_{42}\cdot \nu_{\pa \Omega} d \h^2\\
        &  + \int_{(\pa^* T_1 \union \pa^* T_2 \union \pa^* T_4) \intersect \pa S_3 {\color{black}\intersect \pa \Omega}} c_{34}n_{34}\cdot \nu_{\pa \Omega} \,d\h^2 + \int_{\pa^* T_3 \intersect (\pa S_1 \union \pa S_2 \union \pa S_4) {\color{black}\intersect \pa \Omega}} c_{34}n_{43}\cdot \nu_{\pa \Omega} d \h^2.
    \end{aligned}
\end{equation*}
Then, invoking property \eqref{gencij}, we obtain
\begin{equation}\label{11_final_11}
    \begin{aligned}
        E_0(\mt) \geq E_0(\ms) &+ \sum_{1 \leq i <j \leq 4}\int_{\pa S_i \intersect\pa^*T_j  {\color{black}\intersect \pa \Omega}}c_{ij}n_{ij}\cdot \nu_{\pa \Omega}\,d\h^2 \\
        &+ \sum_{1 \leq i <j \leq 4}\int_{\pa S_j \intersect\pa^*T_i {\color{black}\intersect \pa \Omega}}c_{ij}n_{ji}\cdot\nu_{\pa \Omega}\,d\h^2.
    \end{aligned}
\end{equation}
Next we observe that each of the twelve integrals in \eqref{11_final_11} is non-negative by property \eqref{goodnormals} and \cref{lemma1}. To conclude, from \eqref{11_final_11}, we see that unless
\begin{equation}\label{Peter}
    \sum_{i=1}^4\h^2 \Big[ (\pa^* T_i \Delta \pa S_i) \intersect \pa \Omega \Big] = 0,
\end{equation}
we would find that $E_0(\mt) > E_0(\ms)$, contradicting the minimality of $\mt$.
Consequently, the last line of \eqref{1_final} equals $E_0(\ms)$,
and so we find that

    \begin{align}
       E_0(\ms) \geq E_0(\mt) \geq \sum_{1 \leq i < j \leq 4} \int_{\pa^*T_i \intersect \pa^*T_j} c_{ij} n_{ij} \cdot \nu_{ij}^* \, d\h^2 = E_0(\ms),
       \label{innereq}
    \end{align}
again contradicting the minimality of $\mt$ unless the inequalities in \eqref{innereq} are in fact equalities. We conclude that for all distinct $i, j \in \{1,2,3,4\}$, $n_{ij} = \nu_{ij}^*$, $\h^2$ a.e. on $\pa^* T_i \intersect \pa^* T_j \intersect \Omega$. Hence, with an appeal to \cref{hyperplane-theorem}, we conclude that $\pa^* T_i \intersect \pa^* T_j$ are all planar and by \eqref{Peter}, they agree with $\pa S_i \intersect \pa S_j$ on $\pa\Omega$. Necessarily, $\mt=\ms$ and the proof is complete. 
\end{proof}

\section{Construction of a Suitable Domain }\label{Construct}

In this section we describe one procedure by which a perturbation of $\bbS^2$ can serve as  the boundary of a domain $\Omega$ satisfying \eqref{goodnormals}. It will suffice to describe the construction of one `trough' and one `valley' since the other five troughs and three valleys can be built similarly. 

To this end, let us suppose that, say, $\pa \tilde{S}_1\cap\pa \tilde{S}_2$ lies within the $xz$-plane, where as in Section \ref{defn_of_domain}, the tetradehral cone partitions the unit ball in $\R^3$ into $\{\tilde{S}_1,\ldots,\tilde{S}_4\}$. Suppose further that $\tilde{S}_1$ lies in the region $y>0$ so that the normal $n_{12}$ is given by $(0,-1,0).$ Fixing a small number $r_*$, we will place a  hemispherical `valley' of radius $r_*$ centered at a point on the $z$-axis below $(0,0,1)$. Then away from this valley, the `trough' will lie in a neighborhood of the circular arc $\tilde{\gamma}_{12}$ given by 
\begin{equation}
 \tilde{\gamma}_{12} 
=\{(x,0,\sqrt{1-x^2}):\;2r_*<x<\bar{x}\},\label{g12}
\end{equation}
for some $ \bar{x}>0$. We note that within this trough the set $\pa S_1\cap\pa S_2\cap\pa\Omega$ will agree with $\tilde{\gamma}_{12}$, where $\pa\Omega$ denotes the region perturbed from $\bbS^2$.  The number $\bar{x}$ is unimportant since the focus here is just on constructing one trough, one valley and then connecting them up in a $C^1$ manner while preserving conditions \eqref{goodnormals}.
\vskip.1in
\noindent
\underline{Construction of the trough.} 
For each $x\in(2r_*,\bar{x})$, consider the plane having normal given by the tangent vector to $\tilde{\gamma}_{12}.$ Restricted to each such plane, say $\Pi_x$, the trough will be given by a circular arc of radius $r_*$. To execute this part of the construction, we start with the circular arc in the $yz$-plane given by:
\[
\{(0,y,z):\,z=r_*-\sqrt{r_*^2-y^2},\,\abs{y}<r_*/2\}
\]
Then we apply the rotation matrix
\[
R(x):=\begin{pmatrix}
    \cos\alpha(x)&0&-\sin{\alpha(x)}\\
    0&1&0\\
    \sin{\alpha(x)}&0&\cos{\alpha(x)}
\end{pmatrix}
\]
to this arc, where 
\begin{equation}
\tan{\alpha(x)}=\left(\sqrt{1-x^2}\right)'=-\frac{x}{\sqrt{1-x^2}},\label{tan}
\end{equation}
so that the arc lies in a plane parallel to $\Pi_x$. Finally, we translate this rotated arc so that it sits on $\tilde{\gamma}_{12}$. This results in a trough, given parametrically by a map $F$:
\begin{align}
&F(x,y):=\nonumber\\
&\left(x-(r_*-\sqrt{r_*^2-y^2})\sin{\alpha(x)},\,y\,,
\sqrt{1-x^2}+(r_*-\sqrt{r_*^2-y^2})\cos{\alpha(x)}\right),\nonumber\\
\label{Fdefn}
\end{align}
for $2r_*<x<\bar{x},\;\abs{y}<r_*/2$.

Now for this trough we check the condition \eqref{goodnormals} with $i=1$, $j=2$, $n_{12}=(0,-1,0)$ and 
\begin{equation}\nu_{\pa\Omega}(x,y)=\frac{F_x\times F_y}{\abs{F_x\times F_y}}.\label{cross}\end{equation}
Hence, property \eqref{goodnormals} will hold along the trough provided the second component of $\nu_{\pa\Omega}$ satisfies the conditions
\[
\nu_{\pa\Omega}^{(2)}(x,y)<0\;\mbox{for}\;2r_*<x<\bar{x},\;0<y<r_*/2,
\]
and
\[
\nu_{\pa\Omega}^{(2)}(x,0)=0\;\mbox{for}\;2r_*<x<\bar{x},\;0<y<r_*/2.
\]
A routine calculation using \eqref{tan}, \eqref{Fdefn} and \eqref{cross} yields that
\begin{equation}
\nu_{\pa\Omega}^{(2)}(x,y)=
-y\,\frac{(r_*^2-\sqrt{r_*^2-y^2})\big(1+(1-x^2)^2\big)}{\sqrt{r_*^2-y^2}\,(1-x^2)^{5/2}},\label{nu2}
\end{equation}
and so these two conditions are met.
\vskip.1in\noindent
\underline{Connecting the trough to a valley.} The trough just constructed will ultimately connect in a $C^1$ manner to a lower hemisphere of radius $r_*$, centered at the point
$(0,0,1-2r_*)$, namely 
\begin{equation}
\{(x,y,z):\,x^2+y^2+\big(z-(1-2r_*)\big)^2=r_*^2,\;z\leq 1-2r_*\}.\label{hemi}
\end{equation}
Within this hemisphere, the outer normal is given by
\[
\nu_{\pa\Omega}=\frac{(-x,-y,1-2r_*-z)}{\sqrt{x^2+y^2+\big(z-(1-2r_*)\big)^2}},
\]
so clearly $\nu_{\pa\Omega}^{(2)}<0$ for $y>0$ and $\nu_{\pa\Omega}^{(2)}=0$ for $y=0$ and \eqref{goodnormals} is again  satisfied. 
What remains is to connect the trough to this hemisphere in a $C^1$ manner.  Here it is more convenient to use $y$ and $z$ as parameters, rather than $x$ and $y$. To this end we introduce a smooth interpolating function $g(z)$ that will serve to define $\pa S_1\cap \pa S_2\cap \pa\Omega$ via the curve $(g(z),0,z)$ as a replacement for \eqref{g12}. To ensure $C^1$ contact between this interpolating trough and the previous trough and between the interpolating trough and the hemisphere, we take $g:[1-2r_*,\sqrt{1-4r_*^2}]\to\R$ to be any smooth function satisfying the conditions:
\begin{align*}
    &g(1-2r_*)=r_*,\quad g(\sqrt{1-4r_*^2})=2r_*,\\
 &   g'(1-2r_*)=0,\quad g'(\sqrt{1-4r_*^2})=-\frac{\sqrt{1-4r_*^2}}{2r_*},
\end{align*}
with $g'(z)>0$ and $g''(z)>0$ on $(1-2r_*,\sqrt{1-4r_*^2})$.

Then, analogous to the procedure used to obtain the previous portion of the trough, we obtain this interpolation portion by first rotating the circular arc in the $xy$-plane given by 
\[
x=-r_*+\sqrt{r_*^2-y^2}
\]
through an angle $\beta(z)$ such that $\tan{\beta(z)}=g'(z)$ and then translating it by $(g(z),0,z)$. Again, the resulting surface has cross-sections on planes with normal in the direction $(g'(z),0,1)$ given by circular arcs.
The parametrization for this portion of the trough is 
\begin{align}
&G(y,z):=\nonumber\\
&\left(g(z)-\big(r_*-\sqrt{r_*^2-y^2}\big)\cos{\beta(z)},\, y\,, z+\big( r_*-\sqrt{r_*^2-y^2} \big)\sin{\beta(z)}\right)
,\nonumber
\end{align}
for $1-2r_*\leq z\leq \sqrt{1-4r_*^2}$ and $\abs{y}<r_*/2$. Through a calculation similar to the one yielding the formula \eqref{nu2}, we find that on this portion of $\pa\Omega$ one has
\[
\nu_{\pa\Omega}^{(2)}=-y\frac{\left[ -b(y)\beta'(z)+\cos{\beta(z)}+g'(z)\sin{\beta(z)}\right]}{\sqrt{r_*^2-y^2}},
\]
where $b(y):=-r_*+\sqrt{r_*^2-y^2}$. Since $b(y)<0,\, \beta'(z)>0,\, g'(z)>0$ and $\beta\in [0,\pi/2)$ , we see that the expression in brackets is positive and so the conditions in \eqref{goodnormals} again hold.

We observe that, in particular, for the value $z=1-2r_*$, one has $\beta(1-2r_*)=0$ and consequently,
\[
G(y,1-2r_*)=(\sqrt{r_*^2-y^2},y,1-2r_*)
\]
Comparing this with \eqref{hemi} we see that indeed the trough glues on to the hemisphere precisely along the arc $x^2+y^2=r_*^2,\; z=1-2r_*$  for $\abs{y}<r_*/2.$

\begin{figure}
    \centering
    \includegraphics[width=0.5\linewidth]{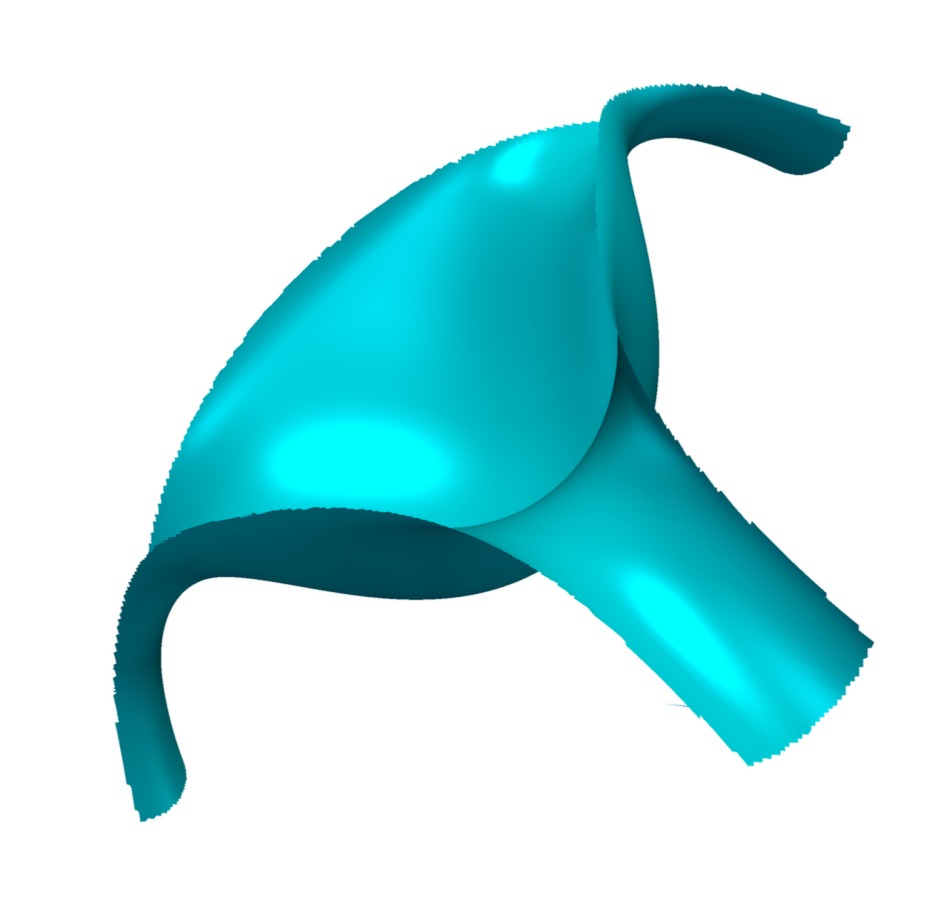}
    \caption{Three troughs meeting a valley.}
    \label{2_troughs_valley}
\end{figure}

This completes the description of how one can construct one trough connecting to a hemispherical valley. This trough is then rotated by $2\pi/3$ in the symmetric case or by the requisite angles in the general case to obtain three troughs terminating at a valley.  The result will resemble \cref{2_troughs_valley}. The whole procedure is repeated to construct all four valleys and six troughs. The rest of the surface constituting $\pa\Omega$ is then just completed in a $C^1$ manner on the complement of the troughs and valleys, forming a domain such as the one in \cref{POI}.  

\vskip.1in\noindent
{\bf Acknowledgments:} We wish to thank Michael Novack for suggesting we pursue a boundary version of the infiltration lemma, which was critical to the proof of Theorem \ref{lemma1}, and for numerous helpful discussions. We thank Matthias Weber {\color{black}and Thanic Nur Samin} for their generous time in producing the figures. The research of P.S. was supported by a Simons Collaboration grant 585520 and an NSF grant DMS 2106516. 
\vskip.1in\noindent
{\bf Conflict of Interest Statement:} The authors have no conflicts of interest to report.
\vskip.1in\noindent
{\bf Data Availability Statement:} There is no external data associated with this research.

\bibliographystyle{alpha} 
\bibliography{refer}

\end{document}